\documentclass[10pt]{article}\usepackage{pb-diagram}

\usepackage{amssymb,latexsym}
\usepackage{amsfonts}
\usepackage{amsmath,amsthm,graphicx}

\newcommand{\bd}{\begin{document}}
\newcommand{\ed}{\end{document}}
\newcommand{\bc}{\begin{center}}
\newcommand{\ec}{\end{center}}
\newcommand{\vs}{\vspace}
\newcommand{\hs}{\hspace}
\newcommand{\bq}{\begin{quote}}
\newcommand{\eq}{\end{quote}}
\newcommand{\mb}{\makebox}
\newcommand{\lt}{\left}
\newcommand{\rt}{\right}
\newcommand{\beqa}{\begin{eqnarray*}}\large
\newcommand{\eeqa}{\end{eqnarray*}}
\newcommand{\beqn}{\begin{eqnarray}}
\newcommand{\eeqn}{\end{eqnarray}}
\newcommand{\bbibl}{\begin{thebibliography}{99}}
\newcommand{\ebibl}{\end{thebibliography}}
\newcommand{\ti}{\times}
\newcommand{\bit}{\begin{itemize}}
\newcommand{\eit}{\end{itemize}}
\newcommand{\ben}{\begin{enumerate}}
\newcommand{\een}{\end{enumerate}}
\newcommand{\lb}{\label}
\newcommand{\hf}{\hspace*{\fill}}
\newcommand{\vf}{\vspace*{\fill}}
\newcommand{\beq}{\begin{equation}}
\newcommand{\eeq}{\end{equation}}
\newcommand{\ba}{\begin{array}}
\newcommand{\ea}{\end{array}}
\newcommand{\del}{\partial}
\newcommand{\bm}[1]{\mb{\boldmath ${#1}$}}
\newcommand{\ot}{\otimes}
\newcommand{\nn}{\nonumber}
\newcommand{\R}{\mb{$I\!\!R$}}
\newcommand{\C}{{\cal C}}
\newcommand{\M}{{\cal M}}
\newcommand{\E}{{\cal E}}
\newcommand{\N}{{\cal N}}
\newcommand{\B}{{\cal B}}
\newcommand{\Y}{{\cal Y}}
\newcommand{\F}{{\cal F}}
\newcommand{\Rc}{{\cal R}}
\newcommand{\A}{{\cal A}}
\renewcommand{\P}{{\cal P}}
\renewcommand{\S}{{\cal S}}
\newcommand{\es}{\emptyset}
\newcommand{\ci}{\subseteq}
\newcommand{\cs}{\supseteq}
\renewcommand{\u}{\cup}
\renewcommand{\i}{\cap}
\newcommand{\bu}{\bigcup}
\newcommand{\bi}{\bigcap}
\newcommand{\la}{\leftarrow}
\newcommand{\ra}{\rightarrow}
\newcommand{\Ra}{\Rightarrow}
\newcommand{\Lra}{\Leftrightarrow}
\newcommand{\lgra}{\longrightarrow}
\newcommand{\Lgra}{\Longrightarrow}
\newcommand{\lglra}{\longleftrightarrow}
\newcommand{\Lglra}{\Longleftrightarrow}
\renewcommand{\a}{\alpha}
\renewcommand{\b}{\beta}
\newcommand{\g}{\gamma}
\newcommand{\G}{\Gamma}
\renewcommand{\d}{\delta}
\newcommand{\D}{\Delta}
\newcommand{\e}{\varepsilon}
\newcommand{\eps}{\epsilon}
\newcommand{\h}{\eta}
\renewcommand{\l}{\lambda}
\newcommand{\m}{\mu}
\newcommand{\n}{\nu}
\newcommand{\p}{\pi}
\newcommand{\s}{\sigma}
\newcommand{\Si}{\Sigma}
\newcommand{\ta}{\tau}
\newcommand{\ph}{\phi}
\newcommand{\Ph}{\Phi}
\renewcommand{\c}{\chi}
\newcommand{\om}{\omega}
\newcommand{\Om}{\Omega}
\newcommand{\tri}{\triangle}
\newcommand{\rec}[1]{\frac{1}{#1}}
\newcommand{\f}{\frac}
\newcommand{\sm}[2]{\sum_{#1}^{#2}}
\newcommand{\ld}{\ldots}
\newcommand{\ov}{\overline}
\newcommand{\ol}[1]{$\bar{\mb{#1}}$}
\newcommand{\un}{\underline}
\newcommand{\iy}{\infty}

\newcommand{\wt}{\widetilde}
\newcommand{\ds}{\displaystyle}
\newcommand{\bdm}{\begin{displaymath}}
\newcommand{\edm}{\end{displaymath}}

\newcommand{\nin}{\not\in}
\newcommand{\bt}{\begin{tabular}}
\newcommand{\et}{\end{tabular}}
\newcommand{\alter}[2]{\lt\{ \ba {ll}#1 \\ #2 \ea \rt.}
\newcommand{\alt}[4]{\lt\{ \ba{ll}#1 & \mb{if \,\,}#2 \\ #3 & \mb{if
               \,\,}#4 \ea \rt.}
\newcommand{\altn}[4]{\lt\{ \ba{rl}#1 & \mb{if \,\,}#2 \\ #3 & \mb{if
               \,\,}#4 \ea \rt.}
\newcommand{\alto}[6]{ \lt\{ \ba{ll}#1 & \mb{if \,\,}#2 \\ #3 & \mb{if
               \,\,} #4 \\ #5 & \mb{if \,\,}#6 \ea \rt.}
\newcommand{\altero}[5]{\mb{$\lt\{ \ba {ll}#1 & \mb{if \,\,}#2 \\ #3 &
               \mb{if \,\,} #4 \\ #5 & \mb{otherwise} \ea \rt.$}}

\newcounter{cnt1}
\newcounter{cnt2}
\newcounter{cnt3}
\newcommand{\blr}{\begin{list}{$($\roman{cnt1}$)$} {\usecounter{cnt1}
        \setlength{\topsep}{0pt} \setlength{\itemsep}{0pt}}}
\newcommand{\bla}{\begin{list}{$($\alph{cnt2}$)$} {\usecounter{cnt2}
        \setlength{\topsep}{0pt} \setlength{\itemsep}{0pt}}}
\newcommand{\bln}{\begin{list}{$($\arabic{cnt3}$)$} {\usecounter{cnt3}
                \setlength{\topsep}{0pt} \setlength{\itemsep}{0pt}}}
\newcommand{\el}{\end{list}}
\newcommand{\no}{\noindent}
\newtheorem{Thm}{Theorem}[section]
\newtheorem{Lem}[Thm]{Lemma}
\newtheorem{Prop}[Thm]{Proposition}
\newtheorem{Def}[Thm]{Definition}
\newtheorem{Exm}[Thm]{Example}
\newtheorem{Rem}[Thm]{Remark}
\newtheorem{Cor}[Thm]{Corollory}
\renewcommand{\baselinestretch}{1}
\newcommand{\ilim}{\mathop{\varprojlim}\limits}
\newcommand{\dlim}{\mathop{\varinjlim}\limits}

\begin{document}
\title{Equivariant Simplicial Cohomology With Local Coefficients and its Classification}
\author{Goutam Mukherjee and Debasis Sen}
\date{}
\maketitle{}
\noindent

\begin{abstract}{We introduce equivariant twisted cohomology of a simplicial set equipped with simplicial action of a discrete group and prove that for suitable twisting function induced from a given equivariant local coefficients, the simplicial version of Bredon-Illman cohomology with local coefficients is isomorphic to equivariant twisted cohomology. The main aim of this paper is to prove a classification theorem for equivariant simplicial cohomology with local coefficients.}
\end{abstract}
\footnote {\bf The second author would like to thank CSIR for its support.}
\section{Introduction} For any abelian group $A$ and natural number $q$, the $q^{th}$ Eilenberg-MacLane simplicial set $K(A,~q)$ represents the $q^{th}$ cohomology group functor with coefficients in $A$ in the sense that for every simplicial set $X$, there is a bijective correspondence \cite{duskin}
$$ H^q(X; A) \cong [X, K(A,~q)].$$

This classification results have been generalized for cohomology with local coefficients in \cite{hir}, \cite{gj}, \cite{bfg}. In this classification the generalized Eilenberg-MacLane complex $L_{\pi}(A,q)$ \cite{git} plays the role of the classifying complex, where the local coefficients system on $X$ is given by an action of $\pi=\pi_1(X)$ on $A$. The complex $L_{\pi}(A,q)$ appears as the total space of a Kan fibration $L_{\pi}(A,q) \longrightarrow \overline{W}\pi$, where $\overline{W}\pi$ denotes the standard $\overline{W}$-construction of $\pi$ \cite{may}. The fibration may be interpreted as an object of the slice category $\mathcal{S}/\overline{W}\pi$, where $\mathcal{S}$ denotes the category of simplicial sets. There is a canonical map $ \theta: X \rightarrow \overline{W}\pi$ and the classification theorem states that the cohomology classes in the $q^{th}$ cohomology with local coefficients of  $X$ correspond bijectively to the vertical homotopy classes of liftings of $\theta.$

The aim of this paper is to prove an equivariant version of the above classifying result. We first describe a simplicial version of Bredon-Illman cohomology with local coefficients \cite{mm}. Next we  define equivariant twisted cohomology of a simplicial set $X$ equipped with simplicial action of a discrete group $G$, where twisting in this equivariant context appears as a natural transformation from the $O_G$-simplicial set $\Phi X$ \cite{elm}, to the category of $O_G$-group complex, $O_G$ being the category of canonical orbits of $G$ \cite{br}. We then prove that for a suitable twisting, induced from a given equivariant local coefficients on $X$, the simplicial version of Bredon-Illman cohomology with local coefficients of $X$ is isomorphic to the equivariant twisted cohomology of $X$.

Finally, we construct a contravariant functor from $O_G$ to the category of Kan complexes which assigns to each object $G/H$ of $O_G$ a generalized Eilenberg-MacLane complex, which will be called the generalized $O_G$-Eilenberg-MacLane complex and use it to prove the classification theorem. We refer \cite{mn}, for  a classification result of Bredon cohomology of simplicial sets equipped with a group action.

The paper is organized as follows. In Section 2, we recall some standard results and fix notations. In Section 3, we introduce the notion of fundamental groupoid and local coefficients system on a simplicial set equipped with simplicial action of a discrete group and define the notion of twisting in the equivariant context. In Section 4, we define the simplicial version of Bredon-Illman cohomology and equivariant twisted cohomology and show that under some connectivity assumptions they are isomorphic. In Section 5, we use a functorial construction of the generalized Eilenberg-MacLane complexes to obtain the generalized $O_G$-Eilenberg-MacLane complex and prove the main classification results.

\section{Preliminaries on $G$-simplicial sets}

In this section we set up our notations, introduce some basic definitions and recall some standard facts.

Let $\Delta$ be the category whose objects are ordered sets
$$[n]=\{0 < 1< \cdots <n\},~n\geq 0,$$ and morphisms are non-decreasing maps $f : [n] \lgra [m].$ There are some distinguished morphisms $d^i:[n-1]\lgra [n], 0\leq i \leq n$, called cofaces and $s^i : [n+1] \lgra [n],~ 0\leq i \leq n$, called codegeneracies, defined as follows:
$$d^i(j) = j,~j<i~~\mbox{and}~~ d^i(j) = j+1,~ j\geq i, ~~(n>0,~  0\leq i \leq n);$$
$$s^i(j) = j,~j \leq i,~~\mbox{and}~~s^i(j) = j-1,~j> i,~~(n \geq 0, ~  0\leq i \leq n).$$
These maps satisfy the standard cosimplicial relations.

A simplicial object $X$ in a category $\mathcal{C}$ is a contravariant functor $X: \Delta \lgra \mathcal{C}.$ Equivalently, a simplicial object is a sequence $\{X_n\}_{n\geq 0}$ of objects of $\mathcal{C}$, together with $\mathcal{C}$-morphisms $ \partial_i : X_n\lgra X_{n-1}$ and $s_i:X_n\lgra X_{n+1},$ $0\leq i \leq n,$
verifying the following simplicial identities:
$$\partial_i \partial_j = \partial_{j-1} \partial_i,~ ~  \partial_i s_j = s_{j-1} \partial_i,~ \mbox{if}~~i<j,$$
$$~~\partial_j s_j =id = \partial_{j+1}s_j,$$
$$\partial_i s_j = s_j \partial_{i-1}, ~~i >j+1;~~ s_is_j = s_{j+1} s_i,~~i\leq j.$$
A simplicial map $f:X\lgra Y$ between two simplicial objects in a category $\mathcal{C}$, is a collection of $\mathcal{C}$-morphisms $f_n : X_n\lgra Y_n,$ $n\geq 0,$  commuting with $\partial_i$ and $s_i$.

In particular, a simplicial set is a simplicial object in the category of sets. Throughout $\mathcal{S}$ will denote the category of simplical sets and simplicial maps.

For any $n$-simplex $x\in X_n$, in a simplicial set $X$, we shall use the notation $\partial_{(i_1, i_2, \cdots, i_r)}x$ to denote the simplex $\partial_{i_1}\partial_{i_2}\cdots \partial_{i_r}x$ obtained by applying the successive face maps $\partial_{i_{r-k}}$ on $x$, where $0\leq i_{r-k}\leq n-k,~0\leq k \leq r-1.$

Recall that the simplicial set $\Delta[n]$, $n\geq 0$, is defined as follows. The set of $q$-simplices is
$$\Delta [n]_q = \{(a_0,a_1,...,a_q) ;~~\mbox{where}~ a_i\in \mathbb{Z},~ 0\leq a_0\leq a_1\leq....\leq a_q\leq n\}.$$
The face and degeneracy maps are defined by $$\partial_i(a_0,...,a_q)=(a_0,..,a_{i-1},a_{i+1},...,a_q),~ s_j(a_0,...,a_q)=(a_0,..,a_i,a_i,..,a_q).$$
Alternatively, the set of $k$-simplices  can be viewed as the contravariant functor
$$\Delta [n]([k]) = ~\mbox{Hom}_{\Delta}([k], [n]),$$
the set of $\Delta$-morphisms from $[k]$ to $[n]$. The only non-degenerate $n$-simplex is $id : [n] \lgra [n]$ and is denoted by $\Delta_n$. In the earlier notation, it is simply, $\Delta_n = (0, 1, \cdots, n).$

It is well known that if $X$ is a  simplicial set, then for any $n$-simplex $x\in X_n$ there is a unique simplicial map $\overline{x} : \Delta[n]\lgra X$ with $\overline{x}(\Delta_n) =x$. Often by an $n$-simplex in a simplicial set $X$ we shall mean either an element $x \in X_n$ or the corresponding simplicial map $\overline{x}.$

We have sinplicial maps $\delta_i:\Delta[n-1]\rightarrow \Delta[n]$ and $\sigma_i:\Delta[n+1]\rightarrow \Delta[n]$ for $0\leq i \leq n$ defined by $\delta_i(\Delta_{n-1})=\partial_i(\Delta_n)$ and $\sigma_i(\Delta_{n+1})=s_i(\Delta_n)$. The boundary subcomplex $\partial \Delta[n]$ of $\Delta[n]$ is defined as the smallest subcomplex of $\Delta[n]$ containing the faces $\partial_i\Delta_n,~~ i=0,1,...,n $. The k-th horn $\Lambda_{k}^{n}$ of $\Delta[n]$ is the subcomplex of $\Delta[n]$ which is generated by all the faces $\partial_i\Delta_n$ except the k-th face $\partial_k\Delta_n$.

\begin{Def}
Let $G$ be a discrete group.  A $G$-simplicial set is a simplicial object in the category of $G$-sets. More precisely, a $G$-simplicial set is a simplicial set $\{X_n ; \partial_i, s_i, 0\leq i \leq n\}_{n\geq 0}$ such that each $X_n$ is a $G$-set and the face maps $\partial_i: X_n\longrightarrow X_{n-1}$ and the degeneracy maps $s_i: X_n\longrightarrow X_{n+1}$ commute with the $G$-action.
A map between a G-simplicial sets is a simplicial map which commutes with the G-action.
\end{Def}

\begin{Def}
A $G$-simplicial set $X$ is called $G$-connected if each fixed point simplicial set $X^H$, $H\subseteq G$, is connected.
\end{Def}
\begin{Def}
Two $G$-maps $f,g:K\rightarrow L$ between two $G$-simplicial sets are $G$-homotopic if there exists a $G$-map $F:K\times \Delta[1]\rightarrow L$ such that
$$F\circ (id\times \delta_1) =f,~~F\circ (id\times \delta_0) =g.$$
The map $F$ is called a $G$-homotopy from $f$ to $g$ and we write $F:f\simeq_G g.$
If $i:K^{\prime}\subseteq K$ is an inclusion of subcomplex and $f,~g$ agree on $K^{\prime}$ then we say that $f$ is $G$-homotopic to $g$ relative to $K^{\prime}$ if there exists a $G$ homotopy $F:f\simeq_G g$ such that $F\circ (i\times id) =\a \circ pr_1,$
where $\alpha=f|_K^{\prime}=g|_K^{\prime}$ and $pr_1 :K^{\prime}\times \Delta[1] \lgra K^{\prime}$ is the projection onto the first factor. In this case we write $F:f\simeq_G g(rel~ K^{\prime}).$
\end{Def}
\begin{Def}
A $G$-simplicial set is a $G$-Kan complex if for every subgroup $H\subseteq G$ the fixed point simplicial set $X^H$ is a Kan complex.
\end{Def}
\begin{Rem}\lb{equi}
Recall (\cite{ag}, \cite{fg}) that the category $G\mathcal{S}$ has a closed model structure \cite{qui}, where the fibrant objects are the $G$-Kan complexes and cofibrant objects are the $G$-simplicial sets. From this it follows that $G$-homotopy on the set of $G$-simplicial maps $K\rightarrow L$ is an equivalence relation, for every $G$-simplicial set $K$. More generally, relative $G$-homotopy is an equivalence relation if the target is a $G$-Kan complex.
\end{Rem}
We consider $G/H \times \Delta[n]$ as a simplicial set where $(G/H \times \Delta[n])_q=G/H \times (\Delta[n])_q$ with face and degeneracy maps as $id\times \partial_i$ and $id \times s_i$. Note that the group $G$ acts on $G/H$  by left translation. With this $G$-action  on the first factor and trivial action on the second factor $G/H \times \Delta[n]$ is a G simplicial set.

A $G$-simplicial map $\sigma:G/H \times \Delta[n]\rightarrow X$ is called an equivariant $n$-simplex of type $H$ in $X$.
\begin{Rem}\lb{correspondence}
We remark that for a $G$-simplicial set $X,$ the set of equivariant $n$-simplices in $X$ is in bijective correspondence with n-simplices of $X^H$. For an equivariant $n$-simplex $\sigma$, the corresponding $n$-simplex is $\sigma^{\prime}=\sigma (eH, \Delta_n).$ The simplicial map $\Delta[n]\lgra X^H,~~\Delta_n \mapsto \sigma^{\prime}$
will be denoted by $\overline{\sigma}.$
\end{Rem}

We shall call $\sigma$  degenerate or non-degenerate according as the $n$-simplex $\sigma^{\prime} \in X^H_n$ is degenerate or non-degenerate.

Recall that the category of canonical orbits, denoted by $O_G,$ is a category whose  objects are cosets $G/H$, as $H$ runs over the all subgroups of $G$. A morphism from $G/H$ to $G/ K$ is a $G$-map. Recall that such a morphism determines and is determined by a subconjugacy relation $g^{-1}Hg\subseteq K$ and is given by $\hat{g}(eH)=gK$. We denote this morphism by $\hat{g}$ \cite{br}.
\begin{Def}
A contravariant functor from $O_G$ to $\mathcal{S}$ is called an $O_G$-simplicial set. A map between $O_G$-simplicial sets is a natural transformation of functors.
\end{Def}
We shall denote the category of $O_G$-simplical sets by $O_G\mathcal{S}.$

\begin{Def}\lb{homotop}
Given two maps $T,S:X\lgra Y$ of $O_G$-simplicial sets, a homotopy $F:T\simeq S$ is defined as follows. Let $X\times \Delta[1]$ be the $O_G$-simplicial set, $G/H\mapsto X(G/H)\times \Delta[1].$ Then a homotopy $F:T\simeq S$ is a map $F:X\times \Delta[1]\lgra Y$ of $O_G$-simplicial sets such that for every object $G/H$ in $O_G$, $F(G/H)$ is a homotopy $T(G/H)\simeq S(G/H)$ of simplicial maps.
\end{Def}
The notion of $O_G$-groups or $O_G$-abelian groups has the obvious meaning replacing $\mathcal{S}$ by $\mathcal{G}rp$ or $\mathcal{A}b.$

If $X$ is a $G$-simplicial set then
$$\underline{X}(G/H):=X^H,~ \underline{X}(\hat{g})=gx,~x\in X^H,~ g^{-1}Hg\subseteq K$$ is an $O_G$-simplicial set. This $O_G$-simplicial set will be denoted by $\Phi X.$

For a $G$-simplicial set $X$, with a $G$-fixed $0$-simplex $v$, we have an $O_G$-group $\underline{\pi}X$ defined as follows. For any subgroup $H$ of $G$,
$$\underline{\pi}X(G/H) := \pi_1(X^H,v)$$ and for a morphism $\hat{g}:G/H\lgra G/K,~~g^{-1}Hg\subseteq K$, $\underline{\pi}X(\hat{g})$ is the homomorphism in fundamental groups induced by the simplicial map $g:X^K\lgra X^H.$
\begin{Def}
An $O_G$-group $\underline{\pi}$ is said to act on an $O_G$-simplicial set (group or abelian group) $\underline{X}$ if for every subgroup $H\subseteq G$, $\underline{\pi}(G/H)$ acts on $\underline{X}(G/H)$ and this action is natural with respect to maps of $O_G.$ Thus if
$$\phi (G/H):\underline{\pi}(G/H)\times \underline{X}(G/H)\lgra \underline{X}(G/H)$$ denotes the action of $\underline{\pi}(G/H)$ on $\underline{X}(G/H)$ then for each subconjugacy relation\\ $g^{-1}Hg\subseteq K,$
$$\phi (G/H)\circ (\underline{\pi}(\hat{g})\times\underline{X}(\widehat{g}))= \underline{X}(\hat{g})\circ \phi(G/K).$$
\end{Def}

\section{Equivariant Local Coefficients and Twisting Function} In this section we introduce the notion of fundamental groupoid of a $G$-simplicial set $X$ and define the notion of (equivariant) local coefficients system on $X$ and discuss a method of generating local coefficients system from an abelian $O_G$-group equipped with an action of the $O_G$-group $\underline{\pi}X.$ At the end of this section we introduce the notion of twisting function in the equivariant context.

We begin with the notion of fundamental groupoid. Recall \cite{gj} that the fundamental groupoid $\pi X$ of a Kan complex $X$ is a category having as objects all $0$-simplexes of $X$ and a morphism $ x\longrightarrow y$ in $\pi X$ is a homotopy class of $1$-simplices $\omega : \Delta [1] \longrightarrow X$ rel $\partial \Delta [1]$ such that $\omega \circ \delta_0 = \overline{y}$, $\omega \circ \delta_1 = \overline{x}$. If $\omega_2$ represents an arrow from $x$ to $y$ and $\omega_0$ represents an arrow from $y$ to $z$,  then their composite $[\omega_0]\circ [\omega_2]$ is represented by $\Omega \circ \delta_1$, where the simplicial map $\Omega:\Delta[2]\lgra X$ corresponds to a $2$-simplex, which is determined by the compatible pair $(\omega^{\prime}_0,~~, \omega^{\prime}_2)$. For a simplicial set $X$ the notion of fundamental groupoid is defined via the geometric realization and the total singular functor.

Recall that an equivariant $n$-simplex of type $H$, $H$ being a subgroup of $G$, is a $G$-simplicial map $\sigma : G/H \times \Delta [n] \longrightarrow X.$ Each such  $\sigma $ corresponds to an $n$-simplex $\sigma^{\prime}\in X^H$ and $\overline{\sigma}:\Delta[n]\lgra X^H$ is the simplicial map given by $\overline{\sigma} (\Delta_n)=\sigma^{\prime}=\sigma (eH, \Delta_n) $. Suppose $x_H$ and $y_K$ are equivariant $0$-simplices of type $H$ and $K$, respectively, and $\hat{g}:G/H \rightarrow G/K$ is a morphism in $O_G$, given by a subconjugacy relation $g^{-1}Hg\subseteq K$, $g\in G,$ so that $\hat{g}(eH)=gK$. Moreover suppose that we have an equivariant $1$-simplex $\phi:G/H \times \Delta[1]\rightarrow X$  of type $H$ such that
$$\phi \circ (id\times \delta_1) = x_H,~~\phi \circ (id\times \delta_0)=y_K\circ (\hat{g}\times id).$$
Then, in particular, $\phi^{\prime}$ is a $1$-simplex in $X^H$  such that $\partial_1\phi^{\prime} = x_H^{\prime}$ and $\partial_0\phi^{\prime} = gy_K^{\prime}$. Observe that the $0$-simplex $gy_K^{\prime}$ in $X^H$ corresponds to the composition
$$G/H \times \Delta [0]\stackrel{\hat{g}\times id}{\rightarrow} G/K \times \Delta [0] \stackrel{y_K}{\longrightarrow} X$$ and $\phi $ is a $G$-homotopy $x_H \simeq_G y_K \circ (\hat{g}\times id)$.

\begin{Def}
Let $X$ be a $G$-Kan complex. The fundamental groupoid $\Pi X$ is a category with objects equivariant $0$-simplices
$$x_{H}:G/ H \times \Delta[0] \rightarrow X$$
of type $H$, as $H$ varies over all subgroups of $G$. Given two objects $x_H$ and  $y_K$ in $\Pi X$, a morphism from $x_H \longrightarrow y_K$ is defined as follows. Consider the set of all pairs $(\hat{g},\phi)$ where $\hat{g}:G/H \rightarrow G/K$ is a morphism in $O_G$, given by a subconjugacy relation $g^{-1}Hg\subseteq K$, $g\in G$ so that $\hat{g}(eH)=gK$ and $\phi:G/H \times \Delta[1]\rightarrow X$ is an equivariant $1$-simplex such that
$$\phi \circ (id\times \delta_1) = x_H,~~\phi \circ (id\times \delta_0)=y_K\circ (\hat{g}\times id).$$

The set of morphisms in $\Pi X$ from $x_H$ to $y_K$ is a quotient of the set of pairs mentioned above by an equivalence relation $`\sim `,$ where $(\hat{g}_{1},\phi_{1})\sim(\hat{g}_{2},\phi_{2})$ if and only if  $g_1=g_2=g$ (say) and there exists a $G$-homotopy
$\Phi : G/H \times \Delta [1] \times \Delta [1] \longrightarrow X$ of $G$-homotopies such that $\Phi : \phi_1 \simeq_G \phi_2$ (rel $G/H \times \partial \Delta [1]$). Since $X$ is a $G$-Kan complex, by Remark \ref{equi}, $\sim$ is an equivalence relation. We denote the equivalence class of $(\hat{g},\phi)$ by $[\hat{g},\phi]$. The set of equivalence classes is the set of morphisms in $\Pi X$ from $x_H$ to $y_K$.

The composition of morphisms in $\Pi X$ is defined as follows. Given two morphisms
$$
\begin{diagram}
\node{x_{H}} \arrow{e,t}{[\hat{g}_{1},\phi_{1}]} \node{y_{K}} \arrow{e,t}{[\hat{g}_{2},\phi_{2}]} \node{z_{L}}
\end{diagram}
$$
their composition $[\hat{g}_2, \phi_2]\circ [\hat{g}_1, \phi_1]$ is $[\widehat{g_1g_2}, \psi]: x_H \longrightarrow z_L$, where the first factor is the composition
\[
\begin{diagram}
\node{G/ H} \arrow[2]{e,t}{\hat{g}_{1}} \node[2]{G/ K} \arrow[2]{e,t}{\hat{g}_{2}} \node[2]{G/L}
\end{diagram}
\]
and $\psi: G/H \times \Delta[1] \longrightarrow X$ is an equivariant $1$-simplex of type $H$  as described below. Let $x$ be a $2$-simplex in the Kan complex $X^H$ determined by the compatible pair of $1$-simplices $(x_0=a_1\phi_2^{\prime}, \hat{x}_1, x_2=\phi_1^{\prime})$ so that $ \partial_0x = g_1\phi_2^{\prime}$ and $\partial_2x= \phi_1^{\prime}$. Then $\psi$ is given by $\psi (eH, \Delta_1) = \partial_1x$.
\end{Def}
For a version of fundamental groupoid of a $G$-space, we refer \cite{mm} and \cite{luck}.

The following lemma shows that the composition is well defined.

\begin{Lem}
The equivalence class of $(\widehat{g_1g_2}, \psi)$ does not depend on the choice of the representatives of $[\hat{g}_1, \phi_1]$ and $[\hat{g}_2, \phi_2]$.
\end{Lem}

\begin{proof}
Suppose $[g_i,\phi_i] = [g_i, \lambda_i],~i=1,2.$ Then there exist $G$-homotopies
$\Theta_i:\phi_{i}\simeq_G \lambda_i$ (rel $G/H\times \partial \Delta [1]$) for $i=1,2.$ Let $y$ be a $2$-simplex in $X^H$ determined by the compatible pair $(y_0=g_1\lambda_2^{\prime}, \hat{y}_1 ,y_2=\lambda_1^{\prime})$ as described above so that $\partial_0y=g_1\lambda_2^{\prime}~ and~ \partial_2y=\lambda_1^{\prime}.$ Let $\xi: G/H \times \Delta[1]\lgra X$ be the equivariant $1$-simplex determined by $\xi (eH, \Delta_1) = \partial_1y$. We need to show that $(\widehat{g_1g_2}, \psi)\sim (\widehat{g_1g_2}, \xi).$  Observe that
$\overline{\Theta}_i:\overline{\phi}_{i}\simeq \overline{\lambda}_i$ (rel $\partial \Delta [1]$) for $i=1,2.$
Now consider the right lifting problem
\[
\begin{diagram}
\node{(\Delta[2]\times \partial\Delta[1])\cup (\Lambda^{2}_{1}\times \Delta[1])} \arrow[3]{e,t}{(\overline{y},\overline{x},\overline{\Theta}_1,\overline{\Theta}_2)}\arrow{s} \node[3]{X^{H}}\arrow{s}\\
\node{\Delta[2]\times \Delta[1]}\arrow[3]{e}\node[3]{*}
\end{diagram}
\]
where in the above diagram, the right vertical arrow is a fibration and the left vertical arrow is an anodyne extensions. Therefore the the above right lifting problem has a solution $\overline{F}:\Delta[2]\times \Delta[1]\rightarrow X^{H}$ and the composition of $\overline{F}$ with
$$\delta_{1}\times id:\Delta[1]\times\Delta[1] \rightarrow \Delta[2]\times\Delta[1]$$
is a homotopy from $\overline{\psi}\simeq \overline{\xi},$ (rel $\partial \Delta[1]$). Let $F: G/H \times \Delta[2] \times \Delta[1] \lgra X$ be the $G$-simplicial map determined by $F(eH,s,t) = \overline{F}(s,t).$ Then the composition
$$G/H\times \Delta[1]\times\Delta[1] \stackrel{id\times \delta_{1}\times id}{\lgra} G/H\times \Delta[2]\times\Delta[1]\stackrel{F}{\rightarrow}X $$ is a $G$-homotopy $\psi \simeq_G \xi$ (rel $G/H\times \partial \Delta [1]).$
As a consequence, $$[\widehat{g_1g_2}, \psi] = [\widehat{g_1g_2}, \xi].$$
\end{proof}
Observe that if $X$ is a $G$-simplicial set then $S|X|$ is a $G$-Kan complex, where for any space $Y$, $SY$  denotes the total singular complex and for any simplicial set $X$, $|X|$ denotes the geometric realization of $X$.
\begin{Def}
For any $G$-simplicial set $X$, we define the fundamental groupoid $\Pi X$ of $X$ by $\Pi X := \Pi S|X|.$
\end{Def}

\begin{Rem}\lb{morphism}
If G is trivial then $\Pi X$ reduces to fundamental groupoid $\pi X$ of a simplicial set X. Again, for a fixed H, the objects $x_H$ together with the morphisms $x_{H}\rightarrow y_{H}$ with identity in the first factor, constitute a subcategory of $\Pi X$ which is precisely the fundamental groupoid $\pi X^H$ of $X^{H}$. Moreover, a morphism $[\hat{g},\phi]$ from $x_H$ to $y_K$, corresponds to the morphism $[\overline{\phi}]$ in the fundamental groupoid $\pi X^H$ of $X^H$ from $x_H^{\prime}$  to $ay_K^{\prime}$, where $\overline{\phi}$ is as in \ref{correspondence}. Suppose $\xi$ is a morphism in $\pi X^H$ from $x$ to $y$ given by a homotopy class $[\overline{\omega}],$ where $\overline{\omega}:\Delta[1]\lgra X^H$ represents the $1$-simplex in $X^H$ from $x$ to $y$. Let $x_H$ and $y_H$ be the objects in $\pi X^H$ defined respectively by $$x_H(eH,\Delta_0)=x,~~y_H(eH,\Delta_0)=y.$$ Then we have a morphism $[id, \omega]:x_H\lgra y_H$ in $\Pi X$, where $\omega(eH,\Delta_1) = \overline{\omega}(\Delta_1).$  We shall denote this morphism corresponding to $\xi$ by $b\xi.$
\end{Rem}

\begin{Def}
An equivariant local coefficient system on a $G$-simplicial set $X$ is a contravariant functor from $\pi X$ to the category Ab of abelian groups.
\end{Def}
\begin{Exm}
Let $X$ be a $G$-simplicial set and n$ >$1. For any object $x_H$ in $\pi X$, define $M(x_H)=\pi_n(X^H,x_H(eH))$ and for any morphism $[\hat{g},\phi]:x_H\rightarrow y_K$ define
$$M([\hat{g},\phi])=([\overline{\phi}])^*\circ \pi_n(g)$$ where $g:X^K\rightarrow X^H$ is the left translation by $g$ and $([\overline{\phi}])^*$ is the isomorphism in fundamental groups of $X^H$ induced by the morphism $[\overline{\phi}]$ from $x_H^{\prime}$ to $y_K^{\prime}.$ Then $M$ is an equivariant local coefficient system on $X$.
\end{Exm}

The following discussion gives a generic example of equivariant local coefficient systems on a $G$-connected $G$-simplicial set $X$ having a $G$-fixed $0$-simplex.

Suppose that $v$ is a $G$-fixed $0$-simplex of $X$ and $M$ is an equivariant local coefficient system on $X$. For any subgroup $H$ of $G$, let $v_{H}$ be the object of type $H$ in $\pi X$ defined by
$$v_{H}:G/H\times \Delta[0]\rightarrow X,$$
$$(eH,\Delta_0)\longmapsto v.$$
Then for any morphism $\widehat{g}:G/H\rightarrow G/ K$ in $O_{G}$ given by a subconjugacy relation $g^{-1}Hg\subseteq K$, we have a morphism  $[\widehat{g},k]:v_{H}\rightarrow v_{K}$ in $\pi X,$ where $k :G/H\times \Delta[1] \lgra X$ is given by $k(eH, \Delta_1)= s_0v$.

Define an abelian $O_G$-group $M_{0}:O_{G}\rightarrow Ab$ by $M_{0}(G/H)=M(v_{H})$ and $M_{0}(\widehat{g})=M[\widehat{g},k]$. The abelian  $O_{G}$ group $M_0$ comes equipped with a natural action of the $O_{G}$ group $\underline{\pi}X$ as described below.

Let $\alpha=[\overline{\phi}]\in \underline{\pi}X(G/H)=\pi_{1}(X^{H},v)$. Then the morphism $[id,\phi]:v_{H}\rightarrow v_{H}$ where $\phi(eH,\Delta_1)= \overline{\phi}(\Delta_1),$ is an equivalence in the category $\pi X$. This yields a group homomorphism
$$b:\pi_{1}(X^{H},v)\rightarrow Aut_{\pi X}(v_{H}), ~\a=[\overline{\phi}] \mapsto b(\a)=[id,\phi].$$ We remark that the composition of the fundamental group $\pi_{1}(X^{H},v)$ coincides with the morphism composition of $\pi X$, contrary to the usual notion of composition in the fundamental group. The composition of the map $b$ with the group homomorphism $Aut_{\pi X}(v_{H})\rightarrow Aut_{Ab}(M(v_{H}))$ which sends $a\in Aut_{\pi X}(v_{H})$ to $[M(a)]^{-1}$ defines the action of $\pi_{1}(X^{H},v)$ on $M_{0}(G/H)$. It is routine to check that this action is natural with respect to morphism of $O_{G}$.

Conversely, an abelian $O_G$-group $M_0$, equipped with an action of the $O_G$-group $\underline{\pi}X,$ defines an equivariant local coefficient system $M$ on $X$ as follows, where $X$ is $G$-connected and $v\in X^G$ a fixed $0$-simplex.

For every object $x_H$ of type $H$, define $M(x_{H})=M_0(G/H).$  To define $M$ on morphisms, we choose a $0$-simplex, say $x$, from each $G$-orbit of $X_{0}$ and an
$1$-simplex $\omega_x$ to fix a morphism $[\overline{\omega_x}]$ from $v$ to $x$ in $\pi X^{G_x}.$  For any other point $y$ in the orbit of x, we fix the $1$-simplex $\omega_y= g\omega_x,$ where $y=gx.$ Clearly, $[\overline{\omega}_{y}]$ is a morphism in $X^{G_y},$ from $v$ to $y.$ Observe that if $x\in X^H$ for a subgroup $H$ of $G$, then $[\overline{\omega}_x]$ is also a morphism in $\pi X^H$ from $v$ to $x$, as $H\subseteq G_x$.

Suppose $x_H \xrightarrow{[\hat{g},\phi]} y_K$ is a morphism in $\pi X.$ Then
by Remark \ref{morphism}, we have a morphism $[\overline{\phi}]$ from $x_H^{\prime}$  to $gy_K^{\prime}$ in $\pi X^H.$  Define $M([\hat{g},\phi])$ as the following composition
$$M_0(G/ K)\xrightarrow{M_0(\hat{g})} M_0(G/H) \xrightarrow{\alpha^{-1}} M_0(G/H),$$
where $\a \in \pi_1(X^H,v)$ is
$$\alpha:=[\overline{\omega}_{gy^{\prime}_K}]^{-1}\circ [\overline{\phi}]\circ [\overline{\omega}_{x^{\prime}_H}],$$
a composition of morphisms in $\pi X^H$ and the second arrow denotes the inverse of the given action of $[\a]$ on $M_0(G/H).$ A straightforward verification shows that $M$ is an equivariant local coefficient system on $X$.

Next we introduce the notion of twisting in the equivariant context. We recall the notion of twisting function on a simplicial set \cite{may}.

\begin{Def}
Let $B$ be a simplicial set and $\Gamma$  a group complex. Then a graded function $$\kappa:B\lgra \gamma,~~\kappa_q:B_q\rightarrow \Gamma_{q-1}$$ is called a twisting function if it satisfies the following identities:
\begin{equation*}
\begin{split}
 \partial_0(\kappa_q(b))&=(\kappa_{q-1}(\partial_0b))^{-1}\kappa_{q-1}(\partial_1 b),~~ b\in B_q\\
\partial_i(\kappa_q(b))&=\kappa_{q-1}(\partial_{i+1}b)~~ i>0\\
s_i(\kappa_q(b))&=\kappa_{q+1}(s_{i+1}b) ~~i\geq 0\\
\kappa_{q+1}(s_0 b)&=e_q,~~e_q \mbox{being the identity of the group }\Gamma_q
\end{split}
\end{equation*}
\end{Def}

\begin{Def}
Let $B$ be any $O_G$-simplicial set and $\Gamma$ be an $O_G$-group complex. Then an $O_G$-twisting function is a natural transformation $\kappa:B\rightarrow \Gamma$ such that $\kappa(G/H)=\{\kappa(G/H)_n\}$ is an ordinary twisting function on $B(G/H)$ for all subgroups $H$ of $G$.
\end{Def}
\begin{Exm}\lb{twist}
Let $X$ be a $G$-connected $G$-simplicial set and $v$ be a $G$-fixed $0$-simplex in $X$. Let $\underline{\pi}X : O_G\lgra \mathcal{G}rp$ be the $O_G$-group as defined before. We regard $\underline{\pi}X$ as an $O_G$-group complex in the trivial way, that is, $\underline{\pi}X(G/H)_n = \underline{\pi}X(G/H)$ for all $n$.
We choose a $0$-simplex $x$ on each $G$-orbit of $X_0$ and a $1$-simplex $\omega_x\in X^{G_x}$ such that $\partial_0\omega_x=x,\partial_1\omega_x=v$. For any other $0$-simplex $y$ on the orbit of $x$ we define $\omega_{y}=g\omega_x$ if $y=gx$. Then it is easy check that this is well defined and $\omega_{y}\in X^{G_y}_1$. For a $0$-simplex $x\in X^H,$ let $\xi_H(x)=[\overline{\omega}_x]$ be the homotopy class of $\overline{\omega_x}: \Delta[1]\lgra X^H$.
Define $$\{\kappa (G/H)_n\}:X^H\rightarrow \pi_1(X^H,v)$$ by $$\kappa (G/H)_n(y):=\xi_H(\partial_{(0,2,\cdots,n)}y)^{-1}\circ[\overline{\partial_{(2,\cdots,n)}y}]\circ\xi_H(\partial_{(1,\cdots,n)}y)$$
where $y\in (X^H)_n$. It is standard that $\kappa (G/H)$ is a twisting function on $X^H.$ We verify that
$$\kappa: \Phi X \lgra \underline{\pi}X, ~~G/H \mapsto \kappa (G/H)$$ is  natural.
Suppose $H$ and $K$ are subgroups such that $g^{-1}Hg\subseteq K.$ Let $z\in X^K_n$.
Then $y=gz\in X^H_n.$ Observe that if $x_1,x_2\in X^K_1$ are $1$-simplexes such that $\overline{x_1}\simeq \overline{x_2},$ as  simplicial maps into $X^K$ then $\overline{y_1}\simeq\overline{y_2}$ as simplicial maps into $X^H$ where $y_i =gx_i,~~i=1,2.$ Thus
\begin{equation*}
 \begin{split}
   \kappa (G/H)_n\circ \Phi X(\hat{g})(z)\\
&= \kappa (G/H)_n(y)\\
&= \xi_H(\partial_{(0,2,\cdots,n)}y)^{-1}\circ[\overline{\partial_{(2,\cdots,n)}y}]\circ\xi_H(\partial_{(1,\cdots,n)}y)\\
&= \xi_H(g\partial_{(0,2,\cdots,n)}z)^{-1}\circ[\overline{g\partial_{(2,\cdots,n)}z}]\circ\xi_H(g\partial_{(1,\cdots,n)}z)\\
&= g\xi_K(\partial_{(0,2,\cdots,n)}z)^{-1}\circ g[\overline{\partial_{(2,\cdots,n)}z}]\circ g\xi_K(\partial_{(1,\cdots,n)}z)\\
&= \underline{\pi}X(\hat{g})\circ \kappa (G/K)_n(z).
  \end{split}
\end{equation*}
Thus $\kappa: \Phi X \lgra \underline{\pi}X$ is an $O_G$-twisting function.
\end{Exm}

\section{Equivariant Simplicial Cohomology with \\ Local Coefficients} In this section, we define a simplicial version of Bredon-Illman cohomology with local coefficients (cf. \cite{mm}) and give an alternative description of this cohomology via the notion of twisting as introduced in the last section.

Let $X$ be a $G$-simplicial set and $M$ an equivariant local coefficients system on $X$. For each equivariant $n$-simplex $\sigma:G/H\times \Delta[n]\rightarrow X,$ we associate an equivariant $0$-simplex $\sigma _{H}:G/H\times\Delta[0]\rightarrow X$ given by
$$\sigma_H= \sigma\circ (id\times \delta_{(1,2,\cdots,n}),$$
where $\delta_{(1,2, \cdots ,n)}$ is the composition
$$\delta_{(1,2, \cdots ,n)}: \Delta[0]\stackrel{\delta_1}{\rightarrow} \Delta[1]\stackrel{\delta_2}{\rightarrow}\cdots \stackrel{\delta_ n}{\rightarrow} \Delta[n].$$
The $j$-th face of $\sigma$ is an equivariant $(n-1)$-simplex of type $H$, denoted by $\sigma^{(j)}$, and is defined by
$$ \sigma^{(j)}= \sigma \circ (id \times \delta_j), ~0\leq j\leq n.$$
\begin{Rem}\lb{initial}
Note that $\sigma^{(j)}_{H}=\sigma_{H}\mbox{ for }j> 0,$ whereas
$$\sigma^{(0)}_H = \sigma \circ (id \times \delta_{(0,2,\cdots, n)}).$$
\end{Rem}

Let $C^{n}_{G}(X;M)$ be the group of all functions $f$ defined on equivariant $n$-simplexes $\sigma:G/H\times\Delta[n]\rightarrow X$ such that $f(\sigma)\in M(\sigma _{H})$ with $f(\sigma)=0,$ if $\sigma$ is degenerate. We have a morphism $\sigma_*=[id,\alpha]$ in $\pi X$ from $\sigma_{H}$ to $\sigma^{(0)}_{H}$ induced by $\sigma$, where $\a : G/H \times \Delta[1] \lgra X$ is given by $\a = \sigma \circ (id \times \delta_{(2,\cdots ,n)}).$ Define a homomorphism $$\delta:C_{G}^{n}(X;M)\rightarrow C_{G}^{n+1}(X;M)$$
$$f\mapsto \delta f$$
where for any equivariant $(n+1)$-simplex $\sigma $ of type $H$,
$$(\delta f)(\sigma)=M(\sigma_{*})(\sigma^{(0)})+\Sigma_{j=1}^{n+1}(-1)^{j}f(\sigma^{(j)}).$$
A routine verification shows that $\delta\circ \delta=0.$ Thus $\{C_{G}^{*}(X;M),\delta \}$ is a cochain complex. We are interested in a subcomplex of this cochain complex as defined below.

Let $\eta:G/H\times \Delta[n]\rightarrow X$ and $\tau:G/ K\times \Delta[n]\rightarrow X$ be two equivariant $n$-simplexes. Suppose there exists a $G$-map $\hat{g} :G/H\lgra G/K,~ g^{-1}Hg\subseteq K,$ such that $\tau \circ (\hat{g}\times id ) = \eta.$ Then $\eta$ and $\tau $ are said to be compatible under $\hat{g}$. Observe that if $\eta$ and $\tau$ are compatible as described above then $\eta$ is degenerate if and only if $\tau$ is degenerate. Moreover notice that in this case, we have a morphism
$[\hat{g},k]:\eta_{H}\rightarrow \tau_{K}$ in $\pi X$, where $k = \eta_H \circ (id \times \sigma_0),$ where $\sigma_0:\Delta[1]\lgra \Delta[0]$ is the map as described in Section 2. Let us denote this induced morphism by $g_*$.

\begin{Def}
We define $S_{G}^{n}(X;M)$ to be the subgroup of $C_{G}^{n}(X;M)$ consisting of all functions f such that if $\eta$ and $\tau$ are equivariant n-simplexes in X which are compatible under $\hat{g} $ then $f(\eta)=M(g_{*})(f(\tau))$.
\end{Def}
\begin{Prop}
If $f\in S_{G}^{n}(X;M) ~then~\delta f \in S_{G}^{n+1}(X;M)$
\end{Prop}
\begin{proof}
Suppose $\eta ,~  \tau$ are equivariant $(n+1)$-simplexes of type $H$ and $K$ respectively, which are compatible under $\hat{g}:G/H\lgra G/K,~g^{-1}Hg\subseteq K.$
Then the faces $\eta^{(j)}$ and $\tau^{(j)}$ are also compatible under $\hat{g}$ for all $j$, $0\leq j\leq n+1.$ Moreover, the induced morphism $g_*:\eta^{(j)}_H\lgra \tau^{(j)}_K$ is the same as the induced morphism $g_*:\eta_H\lgra \tau_K$ for $j\geq 1$ (cf. Remark \ref{initial}) and the compositions
$$\eta_H\stackrel{\eta_*}{\rightarrow}\eta^{(0)}_{H}\stackrel{g_{*}}{\rightarrow}\tau^{(0)}_K~~~\mbox{and}~~~ \eta_H\stackrel{g_*}{\rightarrow}\tau_K\stackrel{\tau_*}{\rightarrow}\tau^{(0)}_K$$
are the same. Thus
\begin{equation*}
 \begin{split}
    M(g_*)(\delta f(\tau))
&= M(g_{*})M(\tau_{*})f(\tau^{(0)})+\Sigma_{j=1}^{n+1}(-1)^{j}M(g_{*})f(\tau^{(j)})\\
&= M(\eta_*)M(g_*)f(\tau^{(0)})+\Sigma_{j=1}^{n+1}(-1)^{j}M(g_{*})f(\tau^{(j)})\\
&= M(\eta_*)f(\eta^{(0)})+\Sigma_{j=1}^{n+1}(-1)^{j}f(\eta^{(j)})= \delta f(\eta).
  \end{split}
\end{equation*}
\end{proof}
Thus we have a cochain complex $S_G(X;M) = \{S^n_G(X;M), \delta\}.$
\begin{Def}
Let X be a G-simplicial set with equivariant local coefficient M on it. Then the    $n$-th Bredon-Illman cohomology of $X$ with local coefficients  $M$ is defined by $$H_{G}^{n}(X;M)=H^{n}(S_{G}(X;M)).$$
\end{Def}

Next, we introduce equivariant twisted cohomology and establish its connection with\\ Bredon-Illman cohomology with local coefficients.

Let $X$ be a $G$-simplicial set having a $G$-fixed $0$-simplex $v$. Suppose $M_0$ is any given abelian $O_G$-group equipped with an action $\phi :\underline{\pi}\times M_0\lgra M_0,$ of an $O_G$-group $\underline{\pi}.$  We regard $\underline{\pi}$ as a trivial $O_G$-group complex. Let $\kappa : \Phi X\lgra \underline{\pi}$ be a given $O_G$-twisting function. We define equivariant twisted cohomology of $X$ with coefficients $M_0$ and twisting $\kappa$ as follows.

We have a cochain complex in the abelian category $\mathcal{C}_G$ defined by $$\underline{C}_n(X):O_G\rightarrow Ab,~~G/H\mapsto C_n(X^H;\mathbb{Z}),$$
where
$C_n(X^H;\mathbb{Z})$ denotes the free
abelian group generated by the non-degenerate $n$-simplexes of $X^H$ and for any morphism in $O_G,$ $\hat{g}:G/H\rightarrow G/K,~~g^{-1}Hg\subseteq K,$ $\underline{C}_n(X)(\hat{g})$ is given by the map
$g_*:C_n(X^K;\mathbb{Z})\rightarrow C_n(X^H;\mathbb{Z}),$ induced by the simplicial map $g:X^K\lgra X^H.$ The boundary $\partial_n:\underline{C}_n(X)\rightarrow \underline{C}_{n-1}(X)$ is a natural
transformation defined by $\partial_n(G/
H):C_n(X^H;\mathbb{Z})\rightarrow C_{n-1}(X^H;\mathbb{Z}),$ where $\partial_n(G/H)$ is the ordinary boundary map of the simplicial set $X^H$. Dualising this chain complex in the abelian
category $\mathcal{C}_G$ we get the cochain complex
$$\{ C^*_G(X;M_0)=Hom_{\mathcal{C}_G}(\underline{C}_*(X),M_0),\delta^n\},$$ which defines the ordinary Bredon cohomology
of the $G$-simplicial set $X$ with coefficients $M_0$. To define the twisted cohomology of
the $G$-simplicial set $X$ we modify the co-boundary maps as follows
$$\delta^n_{\kappa}:C_G^{n}(X;M_0)\rightarrow C_G^{n+1}(X;M_0),~~f\mapsto
\delta^n_{\kappa}f$$
where $$\delta^n_{\kappa}f(G/H):C_{n+1}(X^H;\mathbb{Z})\rightarrow M_0(G/H)$$
is given by
$$\delta^n_{\kappa}f(G/H)(x)=(\kappa(G/H)_{n+1}(x))^{-1}f(G/H)(\partial_0 x)+\Sigma_{i=1}^{n+1}(-1)^i f(G/H)(\partial_i x)$$
for $x\in X^H_{n+1}$. Note that the first term of the right hand side is obtained by the given action $\phi$. We denote the resulting cochain complex by
$C^*_G(X;\kappa, \phi).$
\begin{Def}
The $n^{th}$ equivariant twisted cohomology of $X$ is defined by
$$H^n_G(X;\kappa,\phi):=H_n(C^*_G(X;\kappa, \phi)).$$
\end{Def}
\begin{Rem}\lb{ogcohomology}
We remark that this definition of twisted cohomology makes sense when $\Phi X$ is replaced by any $O_G$-simplicial set.
\end{Rem}

Suppose $X$ is a $G$-connected $G$-simplicial set with a $G$-fixed $0$-simplex $v$. Let $M$ be an equivariant local coefficients system on $X$ and $M_0$ be the associated abelian $O_G$-group equipped with an action of the $O_G$-group $\underline{\pi}X$ as introduced in Section $2$. Let $\kappa$ be the $O_G$-twisting function on as introduced in Example \ref{twist}.

\begin{Thm}\lb{iso}
With the above hypothesis
$$H_G^{n}(X;M)\simeq H_G^n(X;\kappa,\phi)$$
for all $n$.
\end{Thm}
\begin{proof}
Define a cochain map
$$S_G^{*}(X;M)\xrightarrow{\Psi^*} C_G^*(X;\kappa,\phi)$$
as follows. Let $f\in S_G^{n}(X;M)$ and $y\in(X^H)_n$ be non-degenerate. Let $\sigma$ be the unique equivariant $n$-simplex of type $H$ such that $\sigma(eH,\Delta_n)=y$. Then $$\Psi^n (f)(G/H):C_n(X^H)\rightarrow M_0(G/H)$$
is given by
$$\Psi^n (f)(G/H)(y)=M(b\xi_H(\partial_{(1,\cdots,n)}y))f(\sigma)).$$

To check that $\Psi^n (f) \in C_G^n(X;\kappa,\phi)$, suppose $g^{-1}Hg\subseteq K.$
Note that if $z\in X^K_n$ and $y=gz$, then $y\in X^H$. Moreover, if $\sigma$ be as above and $\tau$ denotes the unique equivariant $n$-simplex of type $K$ such that $\tau (eK, \Delta_n)=z,$ then the equivariant $n$-simplexes $\sigma$ and $\tau$ are compatible under $\hat{g}$. As $f\in S^n_G(X;M)$, we must have $f(\sigma)=M(\sigma_H\xrightarrow{g_*}\tau_K)f(\tau)$. Therefore by definition of $\Psi^n$ we have
\begin{equation*}
 \begin{split}
  \Psi^n (f)(G/H)(y)&=M(v_H\xrightarrow{b\xi_H(\partial_{(1,\cdots,n)}y)}\sigma_H)f(\sigma)\\
&=M(v_H\xrightarrow{b\xi_H(\partial_{(1,\cdots,n)}y)}\sigma_H)M(\sigma_H\xrightarrow{g_*}\tau_K)f(\tau).
 \end{split}
\end{equation*}

On the other hand,
$$M_0(\hat{g})\Psi^n (f)(G/K)(z)=M_0(\hat{g})M(v_K\xrightarrow{b\xi_K(\partial_{(1,\cdots,n)}z))}\tau_K)f(\tau).$$
Recall that $M_0(\hat{g}) =M([\hat{g}, k])$, where $k :G/H\times \Delta[1] \lgra X$ is given by $k(eH, \Delta_1)= s_0v$ and note that
$$g_*\circ b\xi_H(\partial_{(1,\cdots,n)}y) = b\xi_K(\partial_{(1,\cdots,n)}z)\circ [\hat{g}, k]$$
as composition of morphisms in $\pi X$. Thus $\Psi^n (f) \in C_G^n(X;\kappa,\phi).$

To check that $\Psi^*$ is a cochain map, let $f\in S_G^{n}(X;M)$, $y\in X^H_{n+1}$ and let $\sigma$ be the equivariant $(n+1)$-simplex of type $H$ such that $\sigma (eH,\Delta_{n+1}) =y.$ Observe that the $i$-th face $\sigma^{(i)}$ is such that $\sigma^{(i)}(eH, \Delta_n) = \partial_iy.$ Thus by the definition of the twisted coboundary we have
\begin{equation*}
  \begin{split}
 &\delta(\Psi^n (f))(G/H)(y)\\
=&\kappa (G/H)(y)^{-1}\Psi^n (f)(G/H)(\partial_0 y)+\Sigma_{i=1}^{n+1}(-1)^i\Psi^n (f)(G/H)(\partial_i y)\\
=&\kappa (G/H)(y)^{-1}M(b\xi_H(\partial_{(1,\cdots,n)}\partial_0y))f(\sigma^{(0)})+\Sigma_{i=1}^{n+1}(-1)^{i}M(b\xi_H(\partial_{(1,\cdots,n)}\partial_iy))f(\sigma^{(i)})\\
=&\kappa (G/H)(y)^{-1}M(b\xi_H(\partial_{(1,\cdots,n)}\partial_0y))f(\sigma^{(0)})+\Sigma_{i=1}^{n+1}(-1)^{i}M(b\xi_H(\partial_{(1,\cdots,n+1)}y))f(\sigma^{(i)}).
\end{split}
\end{equation*}
Note that $\partial_{(1,\cdots,n+1)}y=\partial_{(1,\cdots,n)}\partial_iy$ for $i>0.$

On the other hand,
\begin{equation*}
 \begin{split}
&\Psi^{n+1}(\delta f)(G/H)(y)\\
=& M(v_H\xrightarrow{b\xi_H(\partial_{(1,\cdots,n+1)}y)}\sigma_H)(\delta f)(y)\\
=& M(v_H\xrightarrow{b\xi_H(\partial_{(1,\cdots,n+1)}y)}\sigma_H)\{M(\sigma_*)f(\sigma^{(0)})+\Sigma_{i=1}^{n+1}(-1)^i f(\sigma^{(i)})\}.
 \end{split}
\end{equation*}
Therefore we need to compare the first two terms on the left hand side of the above two expressions. By the definition of the action of $\underline{\pi}X$ on $M_0$ and by Example \ref{twist}, we have
\begin{equation*}
\begin{split}
&\kappa (G/H)(y)^{-1}M(v_H\xrightarrow{b\xi_H(\partial_{(1,\cdots,n)}\partial_0y)}\sigma^{(0)}_{H})\\
=& M(b\xi_H^{-1}(\partial_{(0,2,\cdots,n+1)}y)\circ b[\overline{\partial_{(2,\cdots,n+1)}y}]\circ \xi_H(\partial_{(1,\cdots,n+1)}y))M(b\xi_H(\partial_{(1,\cdots,n)}\partial_0y)\\
=& M(b\xi_H(\partial_{(1,\cdots,n)}\partial_0y)\circ b\xi_H^{-1}(\partial_{(0,2,\cdots,n+1)}y)\circ b[\overline{\partial_{(2,\cdots,n+1)}y}]\circ b\xi_H(\partial_{(1,\cdots,n+1)}y))\\
=& M(b[\overline{\partial_{(2,\cdots,n+1)}y}]\circ b\xi_H(\partial_{(1,\cdots,n+1)}y))\\
=&M(b\xi_H(\partial_{(1,\cdots,n+1)}y))M(\sigma_*).
\end{split}
\end{equation*}
Observe that $\partial_{(0,2,\cdots,n+1)}y = \partial_{(1,2,\cdots,n)}\partial_0y.$
Hence $\Psi^*$ is a cochain map.

Next we define a map
$$\Gamma^*:C_G^*(X;\kappa,\phi)\rightarrow S_G^*(X;M)$$ as follows.
Let $f\in C_G^n(X;\kappa,\phi)$ and $\sigma:G/H\times \Delta[n]\rightarrow X$ be a non-degenerate equivariant $n$-simplex of type $H$. Let $y=\sigma(eH,\Delta_n)$. Define
$$\Gamma^n (f)(\sigma):=M(\sigma_H\xrightarrow{b\xi_H^{-1}(\partial_{(1,\cdots,n)}y)}v_H)f(G/H)(y).$$
To show that $\Gamma^n(f) \in S^n_G(X;M),$ suppose $g^{-1}Hg\subseteq K$, and $\sigma$, $\tau$ are non-degenerate equivariant $n$-simplexes in X of type $H$ and $K$ respectively, such that $\sigma$ and $\tau$ are compatible under $\hat{g}:G/H\lgra G/K.$ Let $z=\tau(eK,\Delta_n)$. Then $y=gz.$ Note that

\begin{equation*}
 \begin{split}
&  M(\sigma_H\xrightarrow{g_*} \tau_K)(\Gamma^n (f)(\tau))\\
=& M(\sigma_H\xrightarrow{g_*} \tau_K)M(\tau_K\xrightarrow{b\xi^{-1}_K(\partial_{(1,\cdots,n)}z)} v_K)(f(G/K)(z))\\
=& M(b\xi^{-1}_K(\partial_{(1,\cdots,n)}z)\circ g_*)(f(G/K)(z))
 \end{split}
\end{equation*}
and
$$\Gamma^n (f)(\sigma)= M(\sigma_H\xrightarrow{b\xi_H(\partial_{(1,\cdots,n)}y)} v_H)f(G/H)(y).$$
But by naturality of $f$ we have $f(G/H)(y)=M_0(\hat{g})f(G/K)(z)$, moreover, $$b\xi_K(\partial_{(1,\cdots,n)}z)\circ [\hat{g}, k] = g_*\circ b\xi_H(\partial_{(1,\cdots,n)}y).$$ Hence $\Gamma^n (f)(\sigma)=M(g_*)\Gamma^n (f)(\tau).$ Thus $\Gamma^n (f)\in S_G^{n}(X;M).$ It is straight forward to verify that $\Psi^*$ and $\Gamma^*$ are inverses to each other. This completes the proof of the theorem.
\end{proof}

\section{Classification} The aim of this last section is to prove a classification theorem for simplicial version of Bredon-Illman cohomology with local coefficients as introduced in Section $4$. We first prove a classification theorem for equivariant twisted cohomology of a $G$-simplicial set, generalizing the corresponding non-equivariant result \cite{hir}. We then use the isomorphism in Theorem \ref{iso} to deduce the desired result.

Recall that twisted cohomology or cohomology with local coefficients of any simplicial set is classified by the so called generalized Eilenberg-MacLane complexes \cite{hir, gj, bfg}. We note that the construction of the standard generalized Eilenberg-MacLane complex is functorial. This motivates us to introduce the notion of $O_G$-generalized Eilenberg-MacLane complex and use it to prove the classification theorem. We first recall the notion of a generalized Eilenberg-MacLane complex \cite{git}.

For any abelian group $A$ and any integer $n>1$, let $K(A,n)$ be the standard Eilenberg-MacLane complex. Recall that $K(A,n)$ can be described as a simplicial abelian group in the following way. Consider a simplicial abelian group $C(A, n)$ with $q$-simplices $C(A,n)_q=C^n(\Delta[q]; A),$ the group of normalized $n$-cochains of the simplicial set $\Delta[q]$. For $\mu \in C^n(\Delta[q];A)$, $\partial_i\mu$ and $s_j\mu$ are defined by
$$\partial_i\mu(\a) =\mu (\delta_i(\a)),~~ s_j\mu(\b) =\mu(\sigma_j(\b))$$
for any $n$-simplex $\a\in \Delta[q-1]$ and for any $n$-simplex $\b\in \Delta[q+1]$ respectively, where $\delta_i:\Delta[q-1]\lgra\Delta[q]$ and $\sigma_j:\Delta[q+1]\lgra \Delta[q]$ are the simplicial maps as defined in Section 2.
We have a simplicial group homomorphism
$$\delta^n: C(A,n)\lgra C(A, n+1)$$
defined as follows. For $c\in C(A,n)_q$, $\delta^n c \in C(A, n+1)_q$ is the usual simplicial coboundary. Then
$$K(A,n)_q = Ker~\delta^n = Z^n(\Delta[q];A)$$ the group of normalized $n$-cocycles.
Suppose $\pi$ is any group which acts on $A$, the action being given by $\phi:\pi\lgra Aut(A).$ Then $\pi$ acts on $K(A,n)$ in the following way. For $a\in \pi$ and $\mu\in K(A,n)_q$, $a\mu :=\phi(a)\circ \mu.$ Let $W\pi$ denote the standard free acyclic complex corresponding to the group $\pi$.
Then a model for the generalized Eilenberg-MacLane complex of the type $(A,n,\phi)$ is
$$L_{\pi}(A,n)= (K(A,n)\times W\pi)/\pi$$
where we quotient out by the diagonal action.

\begin{Rem}\lb{nat}
\begin{enumerate}
\item Let $W\pi\lgra \overline{W}\pi$ be the universal $\pi$-bundle \cite{may} and $\kappa (\pi) :\overline{W}\pi\lgra \pi$ be the canonical twisting function
$$\kappa (\pi)([g_1,g_2,\cdots,g_q]) =g_1.$$ Then recall \cite{thur} that $L_{\pi}(A,n)$ can be viewed as a TCP \cite{may}
$$K(A,n)\times_{\kappa (\pi)}\overline{W}\pi\lgra \overline{W}\pi.$$

\item Suppose $(A,\phi)$ is a $\pi$-module and $(A^{\prime},\phi^{\prime})$ is a $\pi^{\prime}$-module. Moreover, suppose that $\a :\pi\lgra \pi^{\prime}$ is a group homomorphism. View $A^{\prime}$  as a $\pi$-module via $\a$. Then any $\pi$-module homomorphism $f : A\lgra A^{\prime}$ induces a map
$$f_*:K(A,n)\times_{\kappa (\pi)}\overline{W}\pi\lgra K(A^{\prime},n)\times_{\kappa (\pi^{\prime})}\overline{W}\pi^{\prime} $$ in the obvious way.
\end{enumerate}
\end{Rem}

By the Remark \ref{nat}, we have an $O_G$-Kan complex as described below which plays the role of the classifying space in the present context.

Let $\underline{\pi}$ be an $O_G$-group and $M_0$ be an $O_G$-abelian group equipped with an action $\phi:\underline{\pi}\times M_0\lgra M_0$. Then by the above description we have an $O_G$-simplicial differential graded abelian group $\{C(M_0,n)\}$ where $$(C(M_0,n)(G/H))_q= C^n(\Delta[q]; M_0(G/H))$$ for every object $G/H$ of $O_G$ and $C(M_0,n)(\hat{g})$ is induced by $M_0(\hat{g})$ for every morphism $\hat{g} :G/H\lgra G/K$. We denote by $K(M_0,n)= Ker (\delta^n: C(M_0,n) \lgra C(M_0,n+1))$ the corresponding $O_G$-simplicial abelian group. Note that $K(M_0,n)$ is an $O_G$-Eilenberg Maclane complex. Similarly, let  $\overline{W}\underline{\pi}$ denote the $O_G$-Kan complex defined by $\overline{W}\underline{\pi}(G/H) = \overline{W}\underline{\pi}(G/H)$ for every objects $G/H \in O_G$ and with the obvious definition of $\overline{W}\underline{\pi}(\hat{g})$ induced by $\underline{\pi}(\hat{g})$ for a morphism $\hat{g}:G/H\lgra G/K$ in $O_G$. Note that we have an $O_G$-twisting function $\kappa (\underline{\pi}) : \overline{W}\underline{\pi}\lgra \underline{\pi},$ given by $\kappa (\underline{\pi})(G/H)([g_1, g_2, \cdots, g_q]) = g_1.$ By Remark \ref{nat}, we have the following $O_G$-Kan complexes.
$$\chi_{\phi}(M_0,n):= C(M_0,n)\times _{\kappa(\underline{\pi})}\overline{W}\underline{\pi}$$
and $$L_{\phi}(M_0,n):= K(M_0,n)\times _{\kappa(\underline{\pi})}\overline{W}\underline{\pi}.$$
We have natural projections onto the second factor
$$\chi_{\phi}(M_0,n)\xrightarrow{p}\overline{W}\underline{\pi},~~L_{\phi}(M_0,n)\xrightarrow{p}\overline{W}\underline{\pi}$$ and we view these $O_G$-Kan complexes as objects over $\overline{W}\underline{\pi}.$
\begin{Def}
We call the $O_G$-Kan complex $L_{\phi}(M_0,n)$ the generalised $O_G$-Eilenberg\\ MacLane complex.
\end{Def}

We shall need the following Lemma.
\begin{Lem}\lb{eh}
For a subgroup $H$ of $G$ and an integer $q\geq 0$ consider the $G$-simplicial set $G/H\times \Delta[q].$ Let $M_0$ be an abelian $O_G$-group with a given action $\phi$ of an $O_G$-group $\underline{\pi}.$ Let $\kappa: \Phi (G/H\times \Delta[q])\lgra \underline{\pi}$ be a twisting function. Then there is a cochain isomorphism $$E_H^*:C_G^{*}(G/H\times\Delta[q];\kappa,\phi)\rightarrow C^*(\Delta[q]; M_0(G/H))$$ which is natural with respect to morphisms in $O_G$.
\end{Lem}
\begin{proof}
Let $f\in C_G^{n}(G/H\times \Delta[q];\kappa, \phi)$ and $\alpha\in \Delta[q]_n$ be non-degenerate. Suppose $\alpha=(\alpha_0,\cdots,\alpha_n)$ where $0\leq\alpha_0<\alpha_1<\cdots<\alpha_n\leq q$. Note that  $f:\underline{C}_n(G/H\times \Delta[q])\lgra M_0$ is a natural transformation. Define $$E_H^n(f)(\alpha)=\kappa(G/H)(eH,(0,\alpha_0))^{-1}f(G/H)(eH,\alpha).$$
Observe that $(eH, (0, \a_0))$ and $(eH, \a)$ are respectively an $1$-simplex and an $n$-simplex in $(G/H\times \Delta[q])^H$ and the right hand side of the above equality is given by the action of $\underline{\pi}(G/H)$ on $M_0(G/H)$.

To check that $E_H^*$ is a cochain map, let $f\in C_G^{n-1}(G/H\times \Delta[q];\kappa, \phi)$ and $\a = (\alpha_0,\cdots,\alpha_n) \in\Delta[q]_n.$ Then
\begin{equation*}
 \begin{split}
& E^n_H(\delta_{\kappa} f)(\alpha)\\
=& \kappa (G/H)(eH,(0,\alpha_0))^{-1}(\delta_{\kappa} f)(G/H)(eH,\alpha)\\
=& \kappa (G/H)(eH,(0,\alpha_0))^{-1}\{\kappa (G/H)(eH,\alpha)^{-1}
f(G/H)(\partial_0(eH,\alpha))\\
+& \Sigma_{i=1}^n (-1)^i f(G/H)(\partial_i(eH,\alpha))\}.
 \end{split}
\end{equation*}
On the other hand,
\begin{equation*}
 \begin{split}
& \delta(E^{n-1}_H f)(\alpha)\\
=& \Sigma_{i=0}^nE^{n-1}_H(f)(\partial_i\alpha)\\
= & \kappa (G/H)(eH,(0,\alpha_1))^{-1}f(G/H)(eH,\partial_0\alpha)\\
+& \Sigma_{i=1}^n(-1)^i\kappa( G/H)(eH,(0,\alpha_0))^{-1}f(G/H)(eH,\partial_i\alpha).
 \end{split}
\end{equation*}
Note that $\partial_i(eH,\alpha)=(eH,\partial_i\alpha)$. Therefore $E_H^*$ will be a cochain map provided we show that
$$\kappa(G/H)(eH,\alpha)\kappa(G/H)(eH,(0,\alpha_0))=\kappa(G/H)(eH,(0,\alpha_1)).$$
We may assume that $\a_0\neq 0$, for if $\a_0 = 0$, then by the  property of twisted function
$$\kappa (G/H)(eH,(0,\alpha_0))=\kappa (G/H)(s_0(eH,(0)))=e_H,$$
the identity of $\underline{\pi}(G/H)$. Moreover,
\begin{equation*}
 \begin{split}
&\kappa (G/H)(eH,(0,\alpha_1))\\
=& \kappa (G/H)(\partial_{(2,\cdots,n)}(eH,\alpha))\\
= &\partial_{(1,\cdots,n-1)}\kappa (G/H)(eH,\alpha)\\
=& \kappa (G/H)(eH,\alpha).
 \end{split}
\end{equation*}
The last equality holds because all the face maps of the group complex $\underline{\pi}(G/H)$ are identity. So suppose $\a_0 \neq 0$. Now observe that
$\alpha=\partial_0\b$ where $\b=(0,\alpha_0,\cdots,\alpha_n)\in \Delta[q]_{n+1}$. So
$\kappa (G/H)(eH,\alpha)=\kappa (G/H)(\partial_0(eH,\b))$. Furthermore, $$\kappa(G/H)(eH,(0,\alpha_0))=\kappa(G/H)(\partial_{(2,\cdots,n+1)}(eH,\b))=\kappa_H(eH,\b).$$
Therefore
\begin{equation*}
 \begin{split}
 &\kappa(G/H)(eH,\alpha)\kappa (G/H)(eH,(0,\alpha_0))\\
= &\kappa (G/H)(\partial_0(eH,\b))\kappa (G/H)(eH,\b)\\
= &\kappa (G/H)(eH,\partial_1 \b)
 \end{split}
\end{equation*}
Now note that $\partial_1 \b=(0,\alpha_1,\alpha_2,\cdots,\alpha_n).$ As a consequence,
$$\kappa (G/H)(eH,(0,\alpha_1))=\kappa (G/H)(\partial_{(2,\cdots,n)}(eH,\partial_1 \b))=\kappa (G/H)(eH,\partial_1 \b).$$

The inverse
$$(E_H^*)^{-1} : C^*(\Delta[q]; M_0(G/H))\lgra C_G^{*}(G/H\times\Delta[q];\kappa,\phi)$$
is defined as  follows. Suppose $c\in C^n(\Delta[q]; M_0(G/H)).$ Then
$$f=(E_H^*)^{-1}(c):\underline{C}_n(G/H\times \Delta[q])\lgra M_0$$ is given by
$$f(G/K)((\hat{a},\a)) = M_0(\hat{a})(\kappa (G/H)(eH,(0,\a_0))c(\a))$$ for any object $G/K$ in $O_G$ and for any $n$-simplex $(\hat{a}, \a)$ in $(G/H\times \Delta[q])^K,$ where
$$\alpha=(\alpha_0,\cdots,\alpha_n)~~ \mbox{with}~~ 0\leq\alpha_0<\alpha_1<\cdots<\alpha_n\leq q.$$ Observe that
$$\hat{a}\in (G/H)^K = Hom_G(G/K,G/H)=mor_{O_G}(G/K,G/H).$$

To prove the last part, suppose $\hat{g}:G/H\lgra G/K$, $g^{-1}Hg\subseteq K$, is a morphism in $O_G.$ Let $\kappa^{\prime} : \Phi(G/K\times \Delta[q])\lgra \underline{\pi}$ be an $O_G$-twisting function. Let $\kappa=\Phi(\hat{g}\times id)\kappa^{\prime}$ be the twisting function induced by the $G$-map $\hat{g}\times id: G/H\times \Delta[q]\lgra G/K\times \Delta[q].$ Let
$$(\hat{g}\times id)_*:C_G^*(G/K\times \Delta[q];\kappa^{\prime},\phi)\rightarrow C_G^*(G/H\times \Delta[q];\kappa,\phi)$$ be the cochain map induced by $\hat{g}\times id$ and let
$$M_0(\hat{g})_*: C^*(\Delta[q]; M_0(G/K))\lgra C^*(\Delta[q]; M_0(G/H))$$ be the map induced by the coefficient homomorphism $M_0(\hat{g}) : M_0(G/K)\lgra M_0(G/H).$
We need to verify that
$$M_0(\hat{g})_* \circ E_K^* =E_H^*\circ (\hat{g}\times id)^*.$$
Let $f\in C_G^*(G/K\times \Delta[q];\kappa^{\prime},\phi)$ and $\a=(\a_0,\cdots, \a_n)$ be a non-degenerate $n$-simplex in $\Delta[q].$ Then
$$M_0(\hat{g})_n \circ E_K^n(f)(\a)= M_0(\hat{g})(\kappa^{\prime}(G/K)(eK,(0,\alpha_0))^{-1}f(G/K)(eK,\alpha)).$$
On the other hand,
\begin{equation*}
 \begin{split}
E_H^n\circ (\hat{g}\times id)^n(f)(\alpha)\\
= & \kappa (G/H)(eH,(0,\alpha_0))^{-1}((\hat{g}\times id)_n(f))(G/H)(eH,\alpha)\\
= &  \kappa (G/H)(eH,(0,\alpha_0))^{-1} f(G/H)(\hat{g}\times id)(eK,\alpha)\\
= & \kappa^{\prime}(G/H)(gK,(0,\alpha_0))^{-1}M_0(\hat{g})f(G/K)(eK,\alpha)\\
= & M_0(\hat{g})(\kappa^{\prime}(G/K)(eK,(0,\alpha_0))^{-1}f(G/K)(eK,\alpha))
 \end{split}
\end{equation*}
By the naturality of $f$ and the twisting function $\kappa^{\prime}$. Hence the above equality holds.
\end{proof}

Suppose $X$ is a $G$-simplicial set. As before, $M_0$ stands for an abelian $O_G$-group with a given action $\phi$ of an $O_G$-group $\underline{\pi}.$ Let $\kappa : \Phi X\lgra \underline{\pi}$ be a given $O_G$-twisting function and
$$\kappa (\underline{\pi}):\overline{W}\underline{\pi}\lgra \underline{\pi}$$ be the $O_G$-twisting function
$$\kappa (\underline{\pi})(G/H)([g_1, g_2, \cdots, g_q]) = g_1.$$
We have a natural map $\theta (\kappa): \Phi X\lgra \overline{W}\underline{\pi}$ defined as follows.
$$X_{q}^{H}\lgra \overline{W}\underline{\pi}(G/H)_q,$$
$$x\mapsto  [\kappa(G/H)_q(x),\kappa (G/H)_{q-1}(\partial_{0}x),\cdots,\kappa (G/H)_1(\partial_{0}^{q-1}x)].$$
Note that $\kappa (\underline{\pi})\circ \theta(\kappa ) = \kappa.$
Let $(\Phi X,~~\chi_{\phi}(M_0,n))$ denote the set of all liftings of the map $\theta(\kappa)$ under $p$. Clearly $(\Phi X,~~\chi_{\phi}(M_0,n))$ has an abelian group structure induced fibre wise from that of the cochain group. Note that we have coboundary $$C(M_0,n)\times_{\kappa(\underline{\pi})}\overline{W}{\underline{\pi}}\xrightarrow{\delta^n\times_{\kappa(\underline{\pi})}id} C(M_0, n+1)\times_{\kappa(\underline{\pi})}\overline{W}{\underline{\pi}}.$$
If we write $f(G/H)(x)=(c,g)$ where
$$x\in X^H_{q},~c\in C^n(\Delta[q];M_0(G/H)),~~ g\in \overline{W}\underline{\pi}(G/H)_q,$$
then $(\delta^n \times _{\kappa(\underline{\pi})}id)f(G/H)(x)=(\delta^nc,g)$. Thus
$$\{(\Phi X,~~\chi_{\phi}(M_0,n)),\delta^n \times _{\kappa(\underline{\pi})}id\}$$
is a cochain complex.
\begin{Thm}\lb{twco}
There is a cochain  isomorphism
$$\Psi^*:\{(\Phi X,~~\chi_{\phi}(M_0,*))({\delta \times _{\kappa(\underline{\pi})}id})\}\cong \{ C^*_G(X;\kappa,\phi),\delta_{\kappa} \}.$$
\end{Thm}
\begin{proof}
Suppose $f\in (\Phi X,~~\chi_{\phi}(M_0,n)).$ Then $\Psi^nf:\underline{C}_nX\lgra M_0$ is a natural transformation defined as follows. Let $G/H$ be any object in $O_G$ and $x\in X^H_n$. Suppose
$$f(G/H)(x)=(c,g),~c\in C^n(\Delta[n];M_0(G/H)),~ g\in (\overline{W}\underline{\pi}(G/H))_n.$$ Then $\Psi^nf(G/H)(x) = c(\Delta_n).$ The naturality of $\Psi^nf$ follows from the that of $f$. The assignment $f\mapsto \Psi^nf$ defines  the homomorphism $\Psi^n.$

To check that $\Psi^*$ is a cochain map, we compute $\Psi^{n+1}(\delta^n \times _{\kappa(\underline{\pi})}id)f.$ As before for $x\in X^H_{n+1}$, if $f(G/H)(x) = (c,g)$, $c\in C^n(\Delta[n+1]; M_0(G/H))$, $g\in (\overline{W}\underline{\pi}(G/H))_{n+1},$  then $(\delta^n \times _{\kappa(\underline{\pi})}id)f(G/H)(x)=(\delta^nc,g)$. Therefore
\begin{equation*}
 \begin{split}
 &\Psi^{n+1}(\delta \times _{\kappa(\underline{\pi})}id)f(G/H)(x)\\
=&\delta c(\Delta_{n+1})\\
=&\Sigma_{i=0}^{n+1}(-1)^ic(\partial_i\Delta_{n+1}).
\end{split}
\end{equation*}

On the other hand, $\delta_{\kappa}(\Psi^n f)(G/H)(x)$
$$=\kappa (G/H)_{n+1}(x)^{-1}\Psi^n f(G/H)(\partial_0 x)+\Sigma_{i=1}^{n+1}(-1)^i\Psi^n f(G/H)(\partial_i x).$$
Since $f(G/H)$ is simplicial we have
$$f(G/H)(\partial_0x)=\partial_0 f(G/H)(x)=(\kappa(\underline{\pi})(G/H)(g)\partial_0 c,\partial_0g),$$ by the definition of the face map $\partial_0$ in $\chi_{\phi}(M_0,n)(G/H).$ Therefore
$$\Psi^n f(G/H)(\partial_0 x)=\kappa(\underline{\pi})(G/H)(g)\partial_0 c(\Delta_n)$$ Now observe that $$g=p(G/H) f(G/H) (x)=\theta(\kappa)(G/H)(x).$$
As a consequence,
$\kappa(\underline{\pi})(G/H)(g)=\kappa(G/H)_{n+1}(x).$  Thus
$$\kappa (G/H)_{n+1}(x)^{-1}\Psi^n f(G/H)(\partial_0 x)= \partial_0 c(\Delta_n)= c(\delta_0(\Delta_n))= c(\partial_0(\Delta_{n+1})).$$
Similarly, for $i>0$, $\Psi^n f(G/H)(\partial_i x)=\partial_i c(\Delta_n)=c(\delta_i(\Delta_n))= c(\partial_i\Delta_{n+1}).$ Therefore $$\delta_{\kappa}(\Psi^n f)= \Psi^{n+1}(\delta^n \times _{\kappa(\underline{\pi})}id)f.$$
Hence $\Psi$ is a chain map.

Conversely, we define a homomorphism
$$\Gamma^n:C^n_G(X;\kappa,\phi)\rightarrow (\Phi X,~~\chi_{\phi}(M_0,n))$$ in the following way. Let
$$T\in C^n_G(X;\kappa,\phi)=Hom_{\mathcal{C}_G}(\underline{C}_n(X),M_0).$$
To define $\Gamma^nT: \Phi X\lgra \chi_{\phi}(M_0,n),$ note that for any subgroup $H$ of $G$ and $x\in X^H_q,$
$$\Gamma^n T(G/H)(x)\in C^n(\Delta[q];M_0(G/H))\times (\overline{W}\underline{\pi}(G/H))_q$$
with $\theta(\kappa)(G/H)(x)$ as the second component, as $\Gamma^nT$ has to be a lift of $\theta(\kappa).$ To determine the first component of $\Gamma^n T(G/H)(x)$ note that the $G$-simplicial map $\sigma: G/H \times \Delta[q]\lgra X,~~ \sigma(eH, \Delta_q) = x$ induces
$$\sigma^*: C^*_G(X;\kappa, \phi)\lgra C^*_G(G/H\times \Delta[q];\kappa\Phi\sigma,\phi).$$
Using the isomorphism $E^*_H$ of Lemma \ref{eh}, we define
$$\Gamma^n T(G/H)(x) =(E^n_H\sigma^*(T),\theta(\kappa)(G/H)(x)).$$
Suppose $\hat{g}:G/H\lgra G/K,~ g^{-1}Hg\subseteq K$ is any morphism in $O_G$. Let $y\in X_q^K$ and $x=gy.$ Suppose $\tau:G/K\times \Delta[q]\lgra X$ is the $G$-simplicial map with $\tau(eK, \Delta_q) = y.$ Then the $G$-simplicial map $\sigma = \tau\circ(\hat{g}\times id)$ corresponds to $x$. Clearly $(\hat{g}\times id)^*\circ \tau^* =\sigma^*,$ where
$$(\hat{g}\times id)^*: C^*_G(G/K\times \Delta[q]; \kappa\Phi\overline{y},\phi) \lgra C^*_G(G/H\times \Delta[q]; \kappa\Phi\overline{x},\phi)$$ is induced by $\hat{g}\times id.$
This observation along with the naturality of $E^*_H$, imply that $\Gamma^nT$ is a natural transformation. It remains to prove that $\Gamma^*$ is the inverse of $\Psi^*.$

Let $f\in C^n_G(X;\kappa,\phi).$ Then $\Psi^n\Gamma^n(f) = f.$ For if $H\subseteq G,$ $x\in X^H_n$ and $\sigma$ be the equivariant $n$-simplex of type $H$ with $\sigma(eH,\Delta_n) = x$ then
\begin{equation*}
\begin{split}
&\Psi^n\Gamma^n(f)(G/H)(x)\\
&=E^n_H(\sigma^*f) (\Delta_n)\\
&=\{\kappa\Phi\sigma(G/H)(eH,(0,0))\}^{-1}(\sigma^*f)(eH,\Delta_n)\\
&=\{\kappa\Phi\sigma(G/H)(s_0(eH,(0)))\}^{-1}f(G/H)(x)\\
&=e_Hf(G/H)(x)=f(G/H)(x).\\
\end{split}
\end{equation*}
The last two equalities follows from the properties of the twisting function $\kappa\Phi\sigma$. It follows that $\Psi^n\Gamma^n=id.$

Next we prove that $\Gamma^n\Psi^n f=f$ for $f:\Phi X\rightarrow \chi_{\phi}(M_0,n),$ a lift of $\theta(\kappa)$. Let $H\subseteq G$ and $x \in X^H_q$. Let $\sigma :G/H\times \Delta[q]\lgra X$ be the equivariant simplicial map such that $\sigma (eH,\Delta_q)= x.$ Then by the definition of $\Gamma^* $ we have
$$\Gamma^n\Psi^n (f)(G/H)(x)=(E^n_H\sigma^*(\Psi^n f),\theta(\kappa)(G/H)(x)).$$
On the other hand, since $f:\Phi X\rightarrow \chi_{\phi}(M_0,n)$ is a lift of $\theta(\kappa),$ $f(G/H)(x) = (c,u)$, where $u=\theta(\kappa) (G/H) (x)),$ for some cochain
$c \in C^n(\Delta[q];M_0(G/H)).$ We show that $c = E^n_H\sigma^*(\Psi^n f).$
Let $\alpha=(\alpha_0,\cdots,\alpha_n)\in \Delta[q]$ be a non-degenerate $n$-simplex. Then $$\alpha=\partial_{(i_1,i_2,\cdots,i_{q-n})}\Delta_q$$ where $0\leq i_1< i_2<\cdots<i_{q-n}\leq q$ and $$\{\alpha_{0},\cdots\alpha_{n},i_{1},\cdots,i_{q-n}\}=\{0,1,2,\cdots,q\}.$$
Then
\begin{equation*}
\begin{split}
&E^n_H(\sigma^*(\Psi^n f))(\alpha)\\
=&\kappa(G/H)_1 \Phi\sigma(G/H)(eH,(0,\alpha_0))^{-1}\sigma^*(\Psi^n f)(G/H)(eH,\alpha)\\
=&\kappa(G/H)_1 \Phi\sigma(G/H)(eH,(0,\alpha_0))^{-1}\Psi^n f(G/H)(\sigma(eH,\alpha))\\
=&\kappa(G/H)_1 \Phi\sigma (G/H)(eH,(0,\alpha_0))^{-1}\Psi^n
f(G/H)(\partial_{(i_1,i_2,\cdots,i_{q-n})}\sigma(eH,\Delta_q))\\
=&\kappa(G/H)_1 \Phi\sigma (G/H)(eH,(0,\alpha_0))^{-1}(\Psi^n
f(G/H)(\partial_{(i_1,i_2,\cdots,i_{q-n})}x).
\end{split}
\end{equation*}
Suppose $\a_0=0$. Then properties of a twisting function imply
$$\kappa \Phi\sigma (G/H)(eH,(0,\alpha_0))=e_H.$$
Moreover, as $f(G/H)$ is simplicial we have
\begin{equation*}
\begin{split}
&f(G/H)\partial_{(i_1,\cdots,i_{q-n})}(x)\\
=& \partial_{(i_1,\cdots,i_{q-n})}f(G/H)(x)\\
=&\partial_{(i_1\cdots,i_{q-n})}(c,u)\\
=&(\partial_{(i_1\cdots,i_{q- n})}c, \partial_{(i_1,\cdots,i_{q-n})}u).
\end{split}
\end{equation*}
Note that since $\alpha_0=0,$ $i_1$ is greater than zero.
Therefore by the definition of $\Psi^*$
\begin{equation*}
 \begin{split}
&E^n_H(\sigma^*(\Psi^n f))(\alpha)\\
=&\Psi^n f(G/H)(\partial_{(i_1,i_2,\cdots,i_{q-n})}x)\\
=&\partial_{(i_1,i_2,\cdots,i_{q-n})}c(\Delta_n)\\
=& c(\delta_{(i_1,i_2,\cdots,i_{q-n})}\Delta_n)\\
=& c(\alpha)
\end{split}
\end{equation*}
On the other hand if $\alpha_0\neq 0$ then we must have $i_0=0$ and therefore
\begin{equation*}
\begin{split}
&f(G/H)(\partial_{(i_1,\cdots,i_{q-n})}x)\\
=& \partial_{(i_1,\cdots,i_{q-n})}f(G/H)(x)\\
=&\partial_{(0,i_2,\cdots,i_{q-n})}(c,u) \\
=&\partial_0(\partial_{(i_2,\cdots,i_{q-n})}c, \partial_{(i_2,\cdots,i_{q-n})}u)\\
=&(\kappa(\underline{\pi})(G/H)(\partial_{(i_2,\cdots,i_{q-n})}u)\partial_{(0,i_2,\cdots,i_{q-n})}c,\partial_{(0,i_2,\cdots,i_{q-n})}u),
\end{split}
\end{equation*}
by the definition of $\partial_0$ in a twisted product. Thus using the definition of $\Psi^*$ we get
$$\Psi^n f(G/H)(\partial_{(i_1,i_2,\cdots,i_{q-n})}x)=
\kappa(\underline{\pi})(G/H)(\partial_{(i_2,\cdots,i_{q-n})}u)\partial_{(0,i_2,\cdots,i_{q-n})}c(\Delta_n).$$
Now observe that
\begin{equation*}
 \begin{split}
  &\kappa(\underline{\pi})(G/H)(\partial_{(i_2,\cdots,i_{q-n})}u)\\
=&\kappa(\underline{\pi})(G/H)(\partial_{(i_2,\cdots,i_{q-n})}\theta(\kappa)(G/H)(\Phi\sigma)(G/H)(eH,\Delta_q))\\
=&\kappa(\underline{\pi})(G/H)\theta (\kappa)(G/H)\Phi\sigma (G/H)(eH,\partial_{(i_2,\cdots,i_{q-n})}\Delta_q)\\
=&\kappa (G/H)_{n+1}\Phi\sigma (G/H)(eH,(0,\alpha_0,\cdots,\alpha_n))\\
=&\kappa (G/H)_1\Phi\sigma (G/H)(eH,(0,\alpha_0)).
\end{split}
\end{equation*}
The last equality holds because $\Phi\sigma (G/H)$ is a simplicial map,
$$(0,\alpha_0)=\partial_{(2,\cdots,n+1)}(0,\alpha_0,\cdots,\alpha_n)$$ and all the face maps of the group complex $\underline{\pi}(G/H)$ are identity maps.

Therefore
$$E^n_H(\sigma^*(\Psi^nf))(\alpha)=\partial_{(0,i_2,\cdots,i_{q-n})}c(\Delta_n)=c(\alpha).$$
\end{proof}

Consider the $O_G$-twisting function $$\kappa(\underline{\pi})p:\chi_{\phi}(M_0,n)\lgra \underline{\pi}.$$
By Remark \ref{ogcohomology}, the twisted cochain complex $C^*_G(\chi_{\phi}(M_0,n);\kappa(\underline{\pi})p, \phi)$ makes sense. We define a cochain
$$u \in C^n_G(\chi_{\phi}(M_0,n);\kappa(\underline{\pi})p, \phi)= Hom_{\mathcal {C}_G}(\underline{C}_n(C^*_G(\chi_{\phi}(M_0,n)), M_0))$$ as follows. For an object $G/H$ in $O_G$,
$$u(G/H):\underline{C}_n(\chi_{\phi}(M_0,n)(G/H)\rightarrow M_0(G/H)$$ is given by $u(G/H)((c,g))= c(\Delta_n)$ where
$$(c,g)\in (\chi_{\phi}(M_0,n)(G/H)_n=C^n(\Delta[n];M_0(G/H))\times_{\kappa(\underline{\pi})_H}\overline{W}_{\underline{\pi}_H}$$ and $u(\hat{g})$ is induced by $M_0(\hat{g})$ for any morphism $\hat{g}:G/H\lgra G/K,~g^{-1}Hg\subseteq K.$ It is easy to check that $u$ as defined above satisfies the required naturality condition and hence is a cochain.
\begin{Def}
We call the cochain $u \in C^n_G(\chi_{\phi}(M_0,n);\kappa(\underline{\pi})p, \phi)$ the fundamental cochain.
\end{Def}
\begin{Rem}\lb{fund}
Suppose $f\in (\Phi X,~~\chi_{\phi}(M_0,n).$ Then for any object $G/H$ in $O_G$,
$$fG/H):X^H\rightarrow C^{n}(\Delta[-];M_0(G/H))\times _{\kappa(\underline{\pi})_(G/H)}(\overline{W}_{\underline{\pi}}(G/H))$$
induces a cochain map ${f(G/H)^*}$ from the cochain complex
$$C^*(C^n(\Delta[-];M_0(G/H))\times _{\kappa(\underline{\pi})(G/H)}\overline{W}\underline{\pi}(G/H); M_0(G/H))$$ to the cochain complex $ C^*(X^H;M_0(G/H))$ and hence
$$f(G/H)^{*}u(G/H)\in C^n(X^H;M_0(G/H))=Hom(C_n(X^H),M_0(G/H)).$$
Therefore for any $x\in X^H_n,$
$$f(G/H)^*u(G/H)(x)=u(G/H)(f(G/H)(x)) =u(G/H)(c,g)=c(\Delta_{n}).$$
Thus $\Psi^n (f)(G/H)(x)=f(G/H)^*u(G/H)(x)$. Hence $\Psi^n f=f^*(u)$, the pull-back of the fundamental cochain $u$ by $f$.
\end{Rem}
\begin{Cor}\lb{cocy1}
For every $n,$
$$\Gamma^n:C^n_G(X;\kappa,\phi)\rightarrow (\Phi X,~~\chi_{\phi}(M_0,n))$$ restricted to cocycles induces isomorphism
$$Z^n_G(X;\kappa,\phi)\cong (\Phi X,~~L_{\phi}(M_0,n)).$$
\end{Cor}
\begin{Def}
Suppose $f,~g\in (\Phi X,~~L_{\phi}(M_0,n)).$ Then $f$ and $g$ are said to be vertically homotopic if there is a homotopy $F: f\simeq g$ of maps of $O_G$-simplicial sets (cf. Definition \ref{homotop}) such that $p\circ F = \theta (\kappa)\circ pr_1,$ where $pr_1: \Phi X\times \Delta[1]\lgra \Phi X$ is the projection onto the first factor.
\end{Def}
\begin{Prop}\lb{cocy2}
Under the isomorphism
$$Z_G^n(X;\kappa,\phi)\xrightarrow{\Gamma^n}(\Phi X,~L_{\phi}(M_0,n)),$$
$f_0,f_1\in Z_{\phi}^n(X;\kappa)$ are cohomologous if and only if $\Gamma^n f_0,\Gamma^n f_1$ are vertically homotopic.
\end{Prop}
\begin{proof}
Suppose $f_0,f_1\in Z_G^n(X;\kappa,\phi)$ are cohomologous. Then $$f_0=f_1+\delta_{\kappa} h$$ for some $h\in C_G^{n-1}(X;\kappa,\phi)$. Let $\kappa_1$ denote the $O_G$-twisting function obtained by composing $\kappa$ and the projection $\Phi X\times \Delta[1]\lgra \Phi X.$ To show that $\Gamma f_0,\Gamma f_1$ are vertically homotopic, it suffices to find $\gamma \in Z^n _G(X\times \Delta[1];\kappa_1, \phi)$ such that $i_0^*(\gamma)=f_0$  and  $i_1^{*}(\gamma)=f_1,$ where $i_0,i_1:X \rightarrow X\times\Delta[1]$ are two obvious inclusions, because, in that case, the image of $\gamma$ under the isomorphism $$\Gamma :Z_G^{n}(X\times \Delta[1];\kappa_1,\phi)\rightarrow (\Phi(X\times \Delta[1]),~~L_{\phi}(M_0,n))$$ will serve as a vertical homotopy between $\Gamma f_0$ and $\Gamma f_1$.

Let $\gamma_0=pr_1^*f_0\in Z_G^n(X\times \Delta[1];\kappa_1,\phi),$ where $$pr_1^*:C_G^*(X;\kappa,\phi)\rightarrow C_G^*(X\times \Delta[1];\kappa_1,\phi)$$ is the cochain map induced by the projection $X\times \Delta[1]\lgra X.$
Clearly $i_{0}^{*}(\gamma_0)=i_{1}^{*}(\gamma_0)=f_0,$ where $i_{0}^{*},i_{1}^{*}:C_G^*(X\times \Delta[1];\kappa_1, \phi)\rightarrow C_G^*(X;\kappa, \phi)$ are the maps induced by $i_0$ and $i_1$ respectively.
Regard $h\in C_G^{n-1}(X;\kappa,\phi)$ as a cochain defined on $i_1(X)$ and we may extend it to a cochain
$$\beta\in C_G^{n-1}(X\times \Delta[1];\kappa_1,\phi)$$
such that $i_{0}^{*}(\beta)=0,~i_{1}^{*}(\beta)=h$.
Set $\gamma=\gamma_0-\delta\beta$. Observe that
$$i_{0}^{*}(\gamma)= i_{0}^{*}(\gamma_0-(\delta_{\kappa_1}\beta))=f_0-\delta_{\kappa}(i_{0}^{*}\beta)=f_0,$$
and similarly,
$$i_{1}^{*}(\gamma)=f_0-\delta_{\kappa}(i_{1}^{*}\beta)=f_0-\delta_{\kappa}h=f_{1}.$$

Conversely, suppose $\Gamma^nf_0$ and $\Gamma^nf_1$ are vertically homotopic. Then they are homotopic in the sense of Definition \ref{homotop} and so $\Gamma^nf_0(G/H)$ and $\Gamma^nf_1(G/H)$ are simplicially homotopic for any subgroup $H$ of $G$. As a consequence, $$\Gamma^nf_0(G/H)^{*} = \Gamma^nf_1(G/H)^{*}.$$ Therefore by the Remark \ref{fund}, $f_{0}=f_{1}$.
\end{proof}

Recall \cite{dk} that the category $O_G\mathcal{S}$ of $O_G$-simplicial sets is a closed model category in the sense of Quillen \cite{qui}. Moreover, recall that if $C$ is an object of a closed model category $\mathcal{C}$, then the category $\mathcal{C}\downarrow C$, the category of objects over $C$ has a closed model structure, induced from that of $\mathcal{C}$ (cf. page 330 \cite {gj}). In particular, the category $O_G\mathcal{S}\downarrow \overline{W}\underline{\pi}$ of objects over $\overline{W}\underline{\pi} \in O_G\mathcal{S}$ is a closed model category. Consequently, the vertical homotopy of liftings of $\theta (\kappa)$ to $L_{\phi}(M_0,n),$ viewed as abstract homotopy  of  morphisms of $O_G\mathcal{S}\downarrow \overline{W}\underline{\pi},$ is an equivalence relation.

From Corollary \ref{cocy1} and Proposition \ref{cocy2}, we obtain the following result.
\begin{Thm}
Suppose $X$ is a $G$-simplicial set and $\kappa:\Phi X\rightarrow \underline{\pi}$ is a twisting function. Then $$H_G^{n}(X;\kappa,\phi)\cong [\Phi X,~~L_{\phi}(M_0,n)]_{\overline{W}\underline{\pi}}$$ where right hand side denote the vertical homotopy class of lifting of the map $\theta(\kappa).$
\end{Thm}

Suppose $X$ is a $G$-connected $G$-simplicial set with a $G$-fixed $0$-simplex $v$ and assume that $M$ a given equivariant local coefficients system on $X$. Let $M_0$ be the associated abelian $O_G$-group equipped with an action $\phi$ of the $O_G$-group $\underline{\pi}=\underline{\pi}X.$ Let $\kappa$ be the $O_G$-twisting function as in Example \ref{twist}. Then from the above theorem and  Theorem \ref{twco}, we obtain the following result.
\begin{Thm} Under the above hypothesis,
$$H^n_G(X;M) \cong [\Phi X,~~L_{\phi}(M_0,n)]_{\overline{W}\underline{\pi}}~~ \mbox{for all}~~ n.$$
\end{Thm}

{\bf Goutam Mukherjee}\\
Indian Statistical Institute, Kolkata-700108, India.\\
e-mail: goutam@isical.ac.in

{\bf Debasis Sen}\\
Indian Statistical Institute, Kolkata-700108, India.\\
e-mail: dsen\_r@isical.ac.in

\end{document}